\documentclass[11pt,leqno]{amsart}
\usepackage{amsthm,amsfonts,amssymb,amsmath,oldgerm}
\usepackage{epsfig}
\numberwithin{equation}{section}
\usepackage[thinlines]{easybmat}
\usepackage[utf8]{inputenc}
\usepackage[T1]{fontenc}

\newcommand{\uU}{\underline{U}}
\newcommand{\uu}{\underline{u}}
\newcommand{\uz}{\underline{z}}

\renewcommand\d{\partial}
\renewcommand\a{\alpha}
\renewcommand\b{\beta}
\renewcommand\o{\omega}
\def\g{\gamma}

\def\l{\lambda}
\def\eps{\varepsilon }
\def\e{\varepsilon}

\renewcommand\d{\partial}

\renewcommand\d{\partial}
\renewcommand\a{\alpha}

\renewcommand\b{\beta}
\renewcommand\o{\omega}
\newcommand\R{\mathbb R}

\def\g{\gamma}

\def\eps{\varepsilon}
\def\e{\varepsilon}

\def\l{\lambda}

\newcommand\br{\begin{remark}}
\newcommand\er{\end{remark}}
\newcommand\bp{\begin{pmatrix}}
\newcommand\ep{\end{pmatrix}}
\newcommand\be{\begin{equation}}
\newcommand\ee{\end{equation}}
\newcommand\ba{\begin{equation}\begin{aligned}}
\newcommand\ea{\end{aligned}\end{equation}}

\newcommand{\bap}{\begin{app}}
\newcommand{\eap}{\end{app}}
\newcommand{\begs}{\begin{exams}}
\newcommand{\eegs}{\end{exams}}
\newcommand{\beg}{\begin{example}}
\newcommand{\eeg}{\end{exaplem}}
\newcommand{\bpr}{\begin{proposition}}
\newcommand{\epr}{\end{proposition}}
\newcommand{\bt}{\begin{theorem}}
\newcommand{\et}{\end{theorem}}
\newcommand{\bc}{\begin{corollary}}
\newcommand{\ec}{\end{corollary}}
\newcommand{\bl}{\begin{lemma}}
\newcommand{\el}{\end{lemma}}
\newcommand{\bd}{\begin{definition}}
\newcommand{\ed}{\end{definition}}
\newcommand{\brs}{\begin{remarks}}
\newcommand{\ers}{\end{remarks}}

\newtheorem{theo}{Theorem}[section]

\newtheorem{rem}[theo]{Remark}

\newtheorem{exams}[theo]{Examples}

\numberwithin{equation}{section}

\newcommand{\RR}{{\mathbb R}}

\newtheorem{theorem}{Theorem}[section]
\newtheorem{proposition}[theorem]{Proposition}
\newtheorem{corollary}[theorem]{Corollary}
\newtheorem{lemma}[theorem]{Lemma}
\newtheorem{definition}[theorem]{Definition}

\newtheorem{example}[theorem]{Example}
\newtheorem{remark}[theorem]{Remark}

\newcommand\cA{{\cal  A}}

\newcommand\cR{{\cal  R}}
\newcommand\cG{{\cal  G}}

\newcommand\cE{{\cal  E}}

\newcommand\cO{{\cal O}}

\newcommand{\N}{\mathbb{N}}

\newcommand{\ul}{\underline}

\newcommand{\supp}{\text{\rm{supp}}}

\newcommand{\beq}{\begin{equation}}
\newcommand{\eeq}{\end{equation}}

\renewcommand{\a}{\alpha}
\renewcommand{\b}{\beta}
\newcommand{\s}{\sigma}
\renewcommand{\cA}{\mathcal{A}}
\renewcommand{\cE}{\mathcal{E}}
\renewcommand{\cG}{\mathcal{G}}
\renewcommand{\cO}{\mathcal{O}}
\renewcommand{\cR}{\mathcal{R}}
\renewcommand{\L}{\Lambda}

\title[Nonlinear stability and high-frequency damping vs. singularity formation]
{Majda and ZND models for detonation: nonlinear stability vs. formation of singularities }

\author{Paul Blochas}
\address{
Univ Rennes, CNRS, IRMAR - UMR 6625, F-35000 Rennes, France}
\email{{\tt paul.blochas@univ-rennes1.fr}}
\thanks{}

\author{Aric Wheeler}
\address{Indiana University, Bloomington, IN 47405}
\email{awheele@iu.edu }
\thanks{Research of A.W. was partially supported
under NSF grant no. DMS-1700279.}
\thanks{Research of P.B was partially supported by the Institut Universitaire de France and by the French region of Brittany.}

\begin{document}




\begin{abstract}
	For the ZND model, we show also singularity formation on the downstream side for arbitrary exponentially-growing weighted norms. For the Majda model on the other hand, we establish for appropriate such weighted norms a set of energy estimates implying not only non-formation of singularities near waves of arbitrary amplitude but also full asymptotic orbital stability for small-amplitude ones.
\end{abstract}

\maketitle

\section{Introduction}\label{s:intro}
This work is about expanding a new area \cite{DR1,DR2,YZ,BR,FR}:
the inviscid global time-asymptotic stability of piecewise smooth solutions of hyperbolic balance laws: specifically, asking whether smooth and small enough initial perturbations, the solution to \eqref{balance} will remain piecewise smooth, without any other discontinuities appearing.
\be\label{balance}
u_t+ f(u)_x= g(u).
\ee This is quite different from the behavior of conservation laws, for which we have that, generically, other shocks are expected to form. Important physical examples include the Saint-Venant equations \cite{Li2,JNRYZ} and detonation \cite{FD}. \\

The most comprehensive contributions concerning the local existence theory for such discontinuous solutions of balance laws were first obtained by Majda in \cite{Maj2,Maj1} and then by M\'etivier in \cite{Met}. See also the book by Benzoni-Gavage-Serre for a more in-depth exposition. Since then, a number of 1-D results have also been obtained; see \cite{Bress} and references therein. Conditions for finite-time blowup have also been explored, see for example \cite{Met}, with discussions of other cases, such as other types of shocks, especially \cite{BS} and references therein about undercompressive and overcompressive shocks (the second ones being studied in the viscous case, while the first ones need to introduce other conditions to complement the conditions of Rankine-Hugoniot).\\ 

In such a framework, it is possible to ask, given a solution with a discontinuity at a single point (or, in higher dimensions, on a hypersurface) to the equation \eqref{balance}, do smooth perturbations of the wave remain smooth outside of the given jump? Here, the goal is to understand the behavior as $t \rightarrow +\infty$ of perturbations of such waves for two systems, the Majda and ZND models, which are both detonation models that will be presented below. \\ 

A crucial part of the analysis of the Majda model is a high-frequency damping estimate. Damping estimates were initiated in \cite{Z1,Z2} and later expanded in \cite{Z3,JNRZ,JZN,RZ}. The existence of a damping estimate is a property of some systems that allows the control of higher order derivatives of the solution by knowledge of bounds on the lower order ones.\\

In certain contexts such as hyperbolic-parabolic systems like in \cite{Z2} and in the fifth chapter of \cite{FS}, this type of result is obtained through arguments similar to those of Kawashima in \cite{Kaw,KSh}. \\

Damping estimates typically take the form
\be\label{damp}
 \mathcal{E}(v)_t\leq -\theta \mathcal{E}(v) + C\|v\|_{L^2_\alpha}^2,
\ee
where $\mathcal{E}(u)$ is an energy equivalent to $\|u\|_{H^s_\alpha}^2$ and $L^2_\alpha$ and $H^s_\alpha$
are weighted $L^2$ and $H^s$ norms with a weight parameter $\a$. Thus, such an estimate effectively controls the $H^s_\alpha$ norm of the perturbation by the $L^2_\alpha$ norm of the perturbation. \\ 

Links with the high-frequency estimates of the resolvent are discussed in \cite{Z2,RZ}. \\

An important point is that they prevent singularity formation. In the case of conservation laws, it is often expected that shocks will form (see \cite{J,Li1}), which can't happen in presence of such damping estimates. In this case, the first order derivative of the perturbation blows up (more precisely, its $L^{\infty}$ norm blows up), but the $L^p$ norms of the perturbation itself do not blow up. \\

High frequency damping estimates \eqref{damp} are tools to close nonlinear iteration. When proving stability results, they essentially reduce the problem to prove low-regularity estimates that may typically be obtained through a Duhamel formulation. See for instance \cite{YZ}.\\
 
Again, this can be contrasted with the classical results on singularity formation \cite{La,J,Li1,Al,Sp} for the case $g\equiv 0$ of hyperbolic conservation laws. Hence damping estimates depend importantly on properties of $g$. For further discussion, see
\cite{DR1,DR2}. \\	

A systematic treatment of damping estimates has been shown in the context of relaxation systems
\cite{MZ,YZ}.
Here our purpose is to explore limitations of this approach in the physically interesting context of detonation models, which are very similar in structure to and
can be viewed as a degenerate type of relaxation models. We begin by showing that, for this general type of problem, the singularity formation shown by John in \cite{J} extends to this case, that is for perturbations around a wave instead of around some constant. An extensive literature on extensions of this result to more general and geometrically involved situations can be found in the books \cite{Al,C,Sp}. For the Majda model, we will show an asymptotic orbital stability result directly thanks to energy estimates in some weighted space for waves of small amplitude, and damping estimates for the general case. For the ZND model, we show that blow-up will occur for some arbitrarily small initial perturbation and how this prevents the kind of damping estimates that we described above. In particular, unlike the Majda model, one cannot stabilize the shock in the ZND model by using exponentially growing weights. 

\subsection{ZND and inviscid Majda models}\label{s:det}

We focus on two closely related classical models of combustion, an inviscid variation of the model from \cite{M} (see also \cite{F}), which we will call the Majda model, and the reactive Euler equations or Zeldovich-von Neumann-Doering (ZND) model. For simplicity, we will focus on one-step reactions. Both of these models can be written abstractly as
\ba\label{eq:abstractsytem}
U_t+f(U)_x&=k\vec{q}\phi(U)z,\\
z_t&=-k\phi(U)z,
\ea
where $U\in\RR^n$ is comprised of various gas-dynamical properties such as velocity, specific volume, and internal energy, $z\in\RR$ is the mass fraction of unburned gas, $\phi(U)$ is an ``ignition function'', which we will take here to be a rough cutoff depending on the temperature of the gas, $\vec{q}$ corresponds to quantities produced by the reaction, in particular the amount of heat released by the reaction, and $k>0$ corresponds to the reaction rate. Note that in the scalar case $q$ is permitted to have either sign, with $q>0$ corresponding to an exothermic reaction and $q<0$ corresponding to an endothermic reaction. In the Majda model, we take $U,q\in\RR$ to be scalar quantities, with $U$ being a ``lumped variable'' representing features of the density, velocity and temperature of the gas. For the sake of concreteness, we write the Majda model as
\ba\label{eq:Majdamodel}
	U_t+f(U)_x&=kq\phi(U)z,\\
	z_t&=-k\phi(U)z.
\ea
We will write the ZND model in Lagrangian coordinates, where $U$ is now taken to be in $\RR^3$ with $U=(v,u,E)$ for $v$ the specific volume of the gas, $u$ the velocity of the gas, and $E$ the specific gas-dynamical energy (that is, $E=e+\frac{1}{2}u^2$ where $e$ is the specific internal energy). The general system \eqref{eq:abstractsytem} for the ZND model now takes the form given in \cite{Z4,Z5}
\ba\label{eq:ZNDEnergyform}
v_t-u_x&=0,\\
u_t+p_x&=0,\\
E_t+(pu)_x&=qk\phi(T)z,\\
z_t&=-k\phi(T)z.
\ea
To complete the system \eqref{eq:ZNDEnergyform}, one needs to relate the temperature $T$ and pressure $p$ to the variables $(v,u,E)$, or equivalently $(v,u,e)$ for specific internal energy. One common choice to complete the system is to use the ideal gas law to define the pressure and temperature as
\ba\label{eq:idealgaslaw}
p(v,u,e)&=\frac{\Gamma e}{v},\\
T(v,u,e)&=\frac{e}{c},
\ea
where $\Gamma>0$ is the Gruneisen constant and $c$ is the specific heat constant. For our purposes though, the specific forms of $p$ and $T$ are not so important. For other possible choices we refer to \cite{Er1}.\\

A right going detonation wave is a traveling (shock) solution $(U,z)$ of \eqref{eq:Majdamodel} or \eqref{eq:ZNDEnergyform} with speed $\s>0$ satisfying
\be
	\lim_{x\to\pm\infty} (U,z)(x,t)=(U_{\pm},z_{\pm}),
\ee
with $z_+=1$ and $z_-=0$. For the existence of such waves in the Majda model see \cite{Lai} and \cite{M}. For the ZND model, a proof of the existence of such waves is given in \cite{GS,W} for some particular choices of $p$ and $T$. 
Physically, the shock is moving from the totally burned region to the totally unburned region. A long standing question, initiated by Erpenbeck in \cite{Er1} for the ZND model, concerns the stability of these detonation waves. For the ZND model, there are partial stability results such as \cite{Z4,Z5} and works cited therein. In \cite{Z4}, it is shown that ZND detonations are spectrally stable in the weak heat release and high overdrive limits by using techniques from asymptotic ODE theory. The weak heat release limit is $q\to 0$ in \eqref{eq:ZNDEnergyform}. The high overdrive limit concerns a different a choice of $\phi(T)$ in \eqref{eq:ZNDEnergyform} than the one we've made here. It is important to note that our blowup theorem for ZND is specific to the inviscid case. \\

Finally, a closely related model that has been studied in the past is the viscous variation of the ZND model, known as the reactive Navier-Stokes, which may be written abstractly as
\ba
	U_t+f(U)_x&=k\vec{q}\phi(T)z+\e (B(U)U_x)_x,\\
	z_t&=-k\phi(T)+\e (C(U,z)z_x)_x.
\ea
For the reactive Navier-Stokes equations, it is known that spectral stability implies nonlinear (orbital) stability \cite{TZ,Z5}. \\

It is also known that spectral stability of the detonation waves of the reactive Navier-Stokes for all $\e>0$ sufficiently small implies spectral stability of the corresponding ZND detonation wave \cite{Z5}, however, it is still open that nonlinear stability of the detonation wave for all small viscosities implies nonlinear stability of the corresponding ZND detonation wave. Finally, study of the spectral stability in the inviscid limit has been partially done in \cite{LyZ}. \\

Furthermore, there have also been studies of the stability of the wave in the viscous Majda model. In \cite{LyRaTZ,Sz}, it shown that nonlinear stability follows from spectral stability with the help of Green functions methods (the second one being focused on waves of small amplitudes), and the study of the spectral stability has been studied in \cite{JYanZ,JYao,LiY} as well as \cite{JLW} for the low-frequency multi-dimensional variation. Furthemore, the vanishing viscosity limit of this problem on the spectral side has been studied in \cite{RV}. \\

Turning to the inviscid Majda model, an early result on the stability of the detonation wave is \cite{Le}. In that paper, \cite{Le} shows that the weak entropy solution to some Riemann problem converges to the detonation wave as $t\to\infty$. In fact, it is proven by first showing the existence of what the author calls a normal solution to the Riemann problem with initial data $(u_{-\infty},0)$ (where $u_{-\infty}=\lim_{x \rightarrow -\infty}\uu(x)$) on $\R^-$ and $(0,1)$ on $\R^+$. In that case, the solution may be discontinuous at $(t,0)$ for some $t>0$ and hence does not follow from our analysis, as we note that the Riemann data is not in general a small perturbation of our wave, and that our work is centered on the preservation of the smoothness on both sides of the shock. We mention that Levy's result does not require $q$ to be small at the expense of only working with very special initial data. As a final note, it can be noted that a large part of the analysis relies on comparing certain solutions and obtaining monotonicity of $\Theta'$ where $\Theta(t)$ is the position of the shock, while we will rely on energy estimates in our case. We will not, in general, have monotonicity of the derivative of the phase and we will not require the initial perturbation of $\uu$ to have a special sign. There are also results available on the spectral stability of the Majda model. For example, \cite{JYanZ} proves spectral stability of the detonation wave for the inviscid Majda model for piecewise constant ignition functions using Evans function techniques.
\subsection{Local existence theory}
In order to study the asymptotic behavior of solutions to hyperbolic equations in the presence of a shock, we first need to recall the following result on the Cauchy problem. It is an adaptation of the results of chapter 4 from \cite{Met}(specifically Theorems 4.1.5 and 4.1.6). We state a theorem that can be applied to both models studied here. Furthermore, from now on, we will consider solutions in the Lax sense, that is, the one developed by A. Majda and G. Métivier in the references cited before \cite{Maj2,Maj1,Met}. \\

We fix an integer $n \geq 1$, an integer $s \geq 2$, two elements $u_+$ and $u_-$ of $\R^n$. With $b,h : \R^n \rightarrow \R^n$ two smooth functions, we study the equation: $$u_t + (b(u))_x = h(u) \, .$$
We are looking for solutions close to a wave (with a discontinuity) satisfying the equation. More precisely, given $\sigma \in \R$ and a smooth function $\uU : \R_- \rightarrow \R^n$ solving the equation $$(db(\uU(x))-\sigma Id_{\R^n})(\uU'(x))=h(\uU(x)) \, ,$$ with $\uU(0)=:u_- \in \R^n$ and for a given $u_+ \in \R^n$ such that $h(u_+)=0$, we have that $b(u_+)-b(u_-)=\sigma (u_+-u_-)$. We also assume that we have that $\uU$ decays exponentially fast to its limit state, as well as all of its derivatives decay exponentially fast to 0, and that the shock associated to $u_-$ and $u_+$ is stable. We further assume there exists a neighborhood $\mathcal{U}_+$ of $u_+$ and a neighborhood $\mathcal{U}_-$ of $u_-$ such that for $w$ in $\mathcal{U}_+ \cup \mathcal{U}_-$, $\partial_t + db(w)\partial_x$ is constantly hyperbolic (in the $t$ direction). We will be looking for solutions that can be written in the form $t \mapsto (\uU + v(t,\cdot))(\cdot-\phi(t))$ with $v$ in $C^0([0,T),H^s(\R^*))$ with $T \in \R^+ \cup \lbrace + \infty \rbrace$ and $s$ big enough.\\

\textbf{Notation:} As in \cite{Met}, we use the following notation: $CH^s((0,T) \times \R^*)$ where $s$ is a nonnegative integer, $T$ a positive number or $+\infty$ is the subset of $C^0(\overline{[0,T)},H^s(\R^*))$ such that for every $j$ a nonnegative integer with $j \leq s$, we have $u_t \in C^{s-j}(\overline{[0,T)},H^j(\R^*))$ where $\overline{[0,T)}$ is $[0,T]$ if $T<+\infty$ and $\R_+$ otherwise.
\begin{theorem}\label{th1}
In the above framework, there exists $\rho > 0$ such that for every $v_0 \in H^{s+ \,1/2}(\R)^n$ supported away from the shock such that $\|v_0\|_{L^{\infty}} \leq \rho$, there exists $T \in (0,+\infty]$ a unique maximum solution to the equation $u_t+(b(u))_x=h(u)$ in $CH^s((0,T) \times \R^*)$ with phase $\phi$ in $C^{s+1}([0,T))$ with initial data $\uU+v_0$ for $x<0$ and $u_++v_0$ for $x>0$. Furthermore, either $T=+\infty$ or $\limsup_{t \rightarrow T}\|u(t)\|_{L^{\infty}}(\R^*) \geq \rho$ or $\limsup_{t \rightarrow T}\|u_x(t)\|_{L^{\infty}(\R^*)}=+\infty$. \\
Finally, if we have $T \in \R^+$, then a sequence of smooth initial data $(v_n)_n$ that converges (strongly) in $H^s$ to $v \in H^s(\R)$ localized outside of a neighborhood of the point $0$, with $v_n$ giving rise to a solution $u_n$ and $u$ for $v$ with all of them defined on $[0,T]$, we have that $(u_n(t))_n$ converges in $L^2$ to $u(t)$ for every $t \in [0,T]$. 
\end{theorem}
\begin{remark} The two adaptations needed to go from the proof presented in \cite{Met} and the one needed here are: first, we need to add a source term, which barely changes the estimates, and, also to transform the constant state on $\R^-$ to a wave, which can be done by using the finite-speed of propagation of the equation. Furthermore, in (\cite{Met}) the author proves more precise results that are not needed here.
\end{remark}
\subsection{Main results}\label{s:main}
In the positive direction, we have a result of asymptotic orbital stability for the Majda model with weighted norms. In fact, we have the following result (with $H^k_{\eps}(\R^*)$ (where $k \in \N$ and $\eps >0$) being the Sobolev space defined as: $\lbrace v \in H^k(\R^*) \, | \, \forall 0 \leq l \leq k \, , (\partial^l_xv)^2 \exp(\eps |\cdot|) \in L^1 \rbrace$)
\begin{theorem}\label{thm:MajdaStab}
Fix $k > 0$, $f : \R \rightarrow \R$, smooth and $u_0 \in \R^+$ and assume that for some $q_1>0$ and for every $q\in [-q_1,q_1]$, there exists some wave satisfying the requirements of Proposition \ref{exstwave}. Then there exists a $\delta_0 >0$, $q_0 > 0$, $\vartheta >0$, $C>0$ and $\eps >0$ such that for every $q \in [-q_0,q_0]$ and every $(v_0,\zeta_0) \in (H^{\,5/2}(\R) \cap H^2_{\eps}(\R))^2$ supported away from $0$ with $\|v_0\|_{H^2_{\epsilon}}+\|\zeta_0\|_{H^2_{\epsilon}}<\delta_0$, the solution $(u,z)$ to \eqref{eq:Majdamodel} with initial data $(\ul{u}+v_0,\ul{z}+\zeta_0)$ is defined for all $t\in\R_+$. The position of the shock at time $t$, $\psi(t)$, is $C^1$ and for all $t \geq 0$ $$\|u(t,\cdot + \psi(t))-\uu\|_{H^2_{\epsilon}(\R^*)}+\|z(t,\cdot + \psi(t))-\uz\|_{H^2_{\eps}(\R^*)}+|\psi'(t)-\sigma| \leq C (\|v_0\|_{H^2_{\epsilon}}+\|\zeta_0\|_{H^2_{\epsilon}})e^{-\vartheta t}.$$
\end{theorem}
In the negative direction, we have a generalization of the blow up theorem in \cite{J} for perturbations of rapidly decaying shocks. See also \cite{LXY} and references cited therein for related results.
\begin{theorem}\label{thm:generalblowup}
	Let $A(u)$ be strictly hyperbolic with at least one genuinely nonlinear field. Consider a stationary shock solution $\ul{U}$ of
	\be\label{eq:mastereqn}
		U_t+A(U)U_x=0.
	\ee
	Assume that $\ul{U}$ is smooth for $x\not=0$ and the bounds
	\be
		|\d_x^n\ul{U}(x)|\leq c_ne^{-c|x|},
	\ee
	hold for some $c>0$ and for all $n=0,1,2,...$. Then, for all $\theta>0$ small enough, there exists a perturbation $\hat{U}(x,0)$ satisfying
	\begin{enumerate}
		\item $\hat{U}(x,0)$ is compactly supported on an interval $I$ of width one and such that the distance from $I$ to 0 is comparable to $\theta^{-1}$.\\
		\item The $C^2$ norm of $\hat{U}(x,0)$ is $\cO(\theta)$.\\
		\item There exists a $T_0\sim\theta^{-1}$ and $T_*\leq T_0$ such that there is a solution $U$ of \eqref{eq:mastereqn} of the form $U=\ul{U}+\hat{U}$, on the time interval $[0,T_*)$. Moreover, the perturbation $\hat{U}(x,t)$ remains bounded in $L^\infty$ for $t\leq T_*$ but
		\be
			\lim_{t\to T_*^-}||\d_x \hat{U}(\cdot,t)||_{\infty}=\infty.
		\ee
	\end{enumerate}
\end{theorem}

As a corollary, we obtain a blowup result for the ZND model.
\begin{corollary}\label{cor:ZNDBlowup}
	Let $\ul{U}=(\ul{v},\ul{u},\ul{E},\ul{z})$ be a right going Neumann shock of the ZND model \eqref{eq:ZNDEnergyform}. Then, for all $\theta>0$ small enough, there exists a perturbation $\hat{U}(x,0)=(\hat{v},\hat{u},\hat{E},\hat{z})(x,0)$ satisfying
	\begin{enumerate}
		\item $\hat{U}(x,0)$ is supported on an interval $I$ of width one with the distance of $I$ from 0 comparable to $\theta^{-1}$ and $I\subset(-\infty,0)$.\\
		\item $\hat{z}(x,0)=0$.
		\item The $C^2$ norm of $\hat{U}(x,0)$ is $\cO(\theta)$.\\
		\item There exists a $T_0\sim\theta^{-1}$ and $T_*\leq T_0$ such that there is a solution $U$ of \eqref{eq:ZNDEnergyform} of the form $U=\ul{U}+\hat{U}$, on the time interval $[0,T_*)$. Moreover, the perturbation $\hat{U}(x,t)$ remains bounded in $L^\infty$ for $t\leq T_*$ but
		\be
		\lim_{t\to T_*^-}||\d_x \hat{U}(\cdot,t)||_{\infty}=\infty.
		\ee
	\end{enumerate}
	Moreover, one can arrange $\d_x\hat{U}(x,0)$ to be ``maximal'' in an outgoing genuinely nonlinear direction. 
\end{corollary}

We note that Corollary \ref{cor:ZNDBlowup} is not an immediate application of Theorem \ref{thm:generalblowup} due the presence of the reaction terms. The key idea of the proof of Corollary \ref{cor:ZNDBlowup} from Theorem \ref{thm:generalblowup} is to note that $\hat{z}=0$ effectively allows one to take $z=0$ in the ZND model, reducing the ZND model to gas dynamics. We also note that Corollary \ref{cor:ZNDBlowup} also prevents any stability result of a similar form to Theorem \ref{thm:MajdaStab} due to the presence of the outgoing undamped mode.

\subsection{Discussion and open problems}\label{s:disc}
One of the main questions that remains unanswered is the stability of waves for which $\uu$ is not necessarily of small amplitudes. As we have obtained high-frequency estimates and that spectral stability results have been obtained, it is of interest to study the nonlinear stability by using the high-frequency damping estimates and the linear stability. \\

In the recent work \cite{LXY}, the authors perform an in depth study of the blow up of initially small in $L^{\infty}$ but large in $W^{1,\infty}$ initial data to systems of conservation laws, as well as the details of such a blow up. \\

Another main question that is not answered here is whether or not the $W^{1,\infty}$-norm blows up for the ZND model for initial perturbations which are small in the weighted space $H^2_\a$. Here, we only show the instability of the wave or the lack of high-frequency damping as long as the solution remains small. It does not give us finite time blow-up. The main issue is that our adaption of the John arguments requires the characteristics to not interact with the shock on a sufficiently large timescale determined by $\|v_x(0)\|_{\infty}$ for $v$ the perturbation of the shock. The data constructed in the proof of Theorem \ref{thm:generalblowup} has $H^2_\a$ norm of size
\begin{equation*}
	||v||_{H^2_\a}\sim \theta e^{\frac{\a}{\theta}}\gg 1.
\end{equation*}
Rescaling $v$ to $V$ so that $||V||_{H^2_\a}\sim\theta$, then gives us an expected blowup time of size
\begin{equation*}
	T_*\sim\frac{e^{\frac{\a}{\theta}}}{\theta},
\end{equation*}
which is more than long enough for the characteristics emanating from an interval of distance $\theta^{-1}$ to interact with the shock.\\

We note that Corollary \ref{cor:ZNDBlowup} implies there can be no $H^s$ damping for $s\geq 2$ by Sobolev embedding. For the case of conservation laws, one can show that there is no $H^1$-damping either. To see this, we recall from the book \cite{Bress} and the articles \cite{BCP,BLY}, and references cited therein, that initial data with small $BV$-norm have unique solutions in $BV$. Since $H^1$ functions with compact support are $BV$ functions, the $H^1$ solution agrees with the $BV$ solution on the time of existence. On the other hand, it is shown in \cite{Li3,Li4} that (suitably rescaled) solutions of conservation laws converge to linear combinations of $N$-waves in $L^1$ as time increases. These two results can be combined to show that there is $W^{1,p}$-blowup for all $p>1$, in particular the $H^1$ norm blows up as well. Briefly, the observation is the variation in the solution $V(u,I)$, for an interval $I$, is bounded from above by
\be
	V(u,I)\leq |I|^{1-\frac{1}{p}}||u_x||_{L^p(I)}\leq |I|^{1-\frac{1}{p}}||u_x||_{L^p(\RR)},
\ee 
for $|I|$ the width of the interval. However, because the $N$-wave has a discontinuity, eventually there is a point $(x_*,t_*)$ where the solution has variation bounded \emph{away} from zero for any interval $I$ containing $x_*$. This then forces the $L^p$ norm of $u_x$ to be infinite at $t_*$ for $p>1$. Interestingly, this argument does not imply that the $W^{1,1}$-norm blows up, and indeed the equality
\be
	||u_x||_{L^1}=||u||_{BV},
\ee
for smooth solutions $u$ seems to imply that the $W^{1,1}$-norm remains finite up to the formation of the shock. \\

It would be interesting to adapt these results to the case with exponentially small perturbations in the coefficients. Another interesting question is how much longer does the $W^{1,p}$ norm take to blowup? It easy to check, and part of the construction, that for compactly supported initial data, as long as the $C^1$ norm of the solution remains finite, the $W^{1,p}$ norm also remains finite. On the other hand, the $C^1$ norm blowing up does not a priori imply that any other $W^{1,p}$ norm blows up. For scalar conservation laws, motivated by \cite{GX,YAJOUT} \cite{CG} provide an asymptotic expansion of the solution in the vanishing viscosity limit in the time period just before shock formation.\\

The key difference in our results between the Majda and ZND models is that the ZND model has an undamped outgoing acoustic mode, whereas the outgoing mode in the Majda model is damped. This suggests that genuinely nonlinear outgoing undamped characteristics play a key role in the formation of singularities in spaces with weighted norms, as incoming characteristics can be handled by having exponential weights which ``trap'' the perturbations near the shock which would then prevent blowup since incoming signals would interact with the shock before they have a chance to blowup. It would be interesting to see if this mechanism is present in systems of hyperbolic balance laws where one assumes that damping is absent at the equilibrium endstate only for some characteristic direction which is outgoing.\\

As in \cite{BR}, one can ask if there is a uniform stability result in the inviscid limit of the viscous Majda model. More precisely, for the viscous Majda model
\ba\label{eq:visMajda}
	U_t+f(U)_x&=qk\phi(U)z+\nu U_{xx},\\
	z_t&=-k\phi(U)z+\nu z_{xx},
\ea
can one show that there is a one parameter family of shocks $(\ul{U}^\nu,\ul{z}^\nu)$ such that $(\ul{U}^\nu,\ul{z}^\nu)$ converges to the detonation wave of the inviscid Majda model as $\nu\to0$, there is a $\vartheta^\nu$ such that the decay estimate in Theorem \ref{thm:MajdaStab} holds for the viscous shock, and such that $\vartheta^\nu$ can be chosen independently of $\nu$ provided $\nu$ is sufficiently small?\\

The final question we ask here is what does stability, in say the class $BV$ of solutions, look like in the cases where the damping estimates fail? Could one have algebraic decay of small perturbations in time? We note that, at the linear level, the Hille-Yosida theorem prevents exponential decay in any exponentially weighted norm, and, as in \cite{MEDP} the theorem of Datko-Pazy prevents any decay at the linear level without a change of topology between the solution and the initial data.\\


\medskip
{\bf Acknowledgement} We thank our thesis advisors L. Miguel Rodrigues and Kevin Zumbrun for suggesting this problem, and for their guidance and discussions.

\section{Stability for Majda's model}\label{s:majda}
In this section we establish the positive result Theorem \ref{thm:MajdaStab}.
\subsection{Existence of the wave}
Let $u_i>0$, and $\phi$ be defined as \begin{align*}
\begin{cases}
\phi(u)= 1$ if $u>u_i, \\
\phi(u)= 0 $ if $u\leq u_i,
\end{cases}
\end{align*} and $f : \R \rightarrow \R$, a smooth function. \\
Let $k>0$, $q>0$. \\
We consider the following system
\begin{align*}
\begin{cases}
u_t =kq \phi(u)z-(f(u))_x, \\
z_t=-kz \phi(u).
\end{cases}
\end{align*}
First, we are interested in traveling waves solutions $(t,x) \mapsto (\uu_0(x-\sigma t),\uz_0(x-\sigma t))$, with a jump (chosen initially at 0), smooth on $\R_-$, $z_0(0^-)=1$ and, furthermore, with $\uu_0$ smooth on $\R_+$ with $\inf_{\R^-}\uu_0>u_i>\sup_{\R_+}\uu_0$, and $\uu_0(0^+)=0$. \\

As said before, the existence of such waves has already been proven in the literature, but we will recall the proof as it gives sharp bounds on the decay rate of the wave on $\R^-$. \\
\begin{proposition}\label{exstwave} \begin{itemize}
\item If a solution to this problem exists, then $\sigma = \dfrac{f(\uu_0(0^-))-f(0)}{u_0}$, and $z_0(x)= \exp \left( \dfrac{kx}{\sigma} \right) $ if $x<0$, $z_0(x)=1$ (this choice can be changed up to changing the constant $q$ by some multiplicative constant) if $x>0$. We also have that $\uu_0$ solves $(f(\uu_0)-\sigma \uu_0)'(x)=z_0(x)$ on $\R_-$ and $\uu_0(x)=0$ if $x>0$.
\item Furthermore, let $u_0 \in (u_i,+\infty)$ and $\sigma := \dfrac{f(u_0)-f(0)}{u_0}$. \\ We denote by $\uu_{-}$ the maximal (smooth) solution on $\R_-$ to \begin{align*}
\uu_{-}'(x)&=\dfrac{kq \, \exp \left( \dfrac{kx}{\sigma} \right) } {f'(\uu_{-})-\sigma}, \\
\uu_-(0)&=u_0.
\end{align*}
Then $(\uu,\uz)$ defines a traveling wave solution to the previous problem such that $\inf_{\R_-}f'(\uu_{-}(x))>\sigma$ if and only if the equation $f(u)=\sigma u-q\sigma$ has a solution in $u \in (u_i,u_0)$ and the biggest such $u$ (denoted $u_{-\infty}$) satisfies $\inf f'([u_{-\infty},u_0]) > \sigma$ and $\sigma>f'(0)$. 
\end{itemize}
\end{proposition}
\begin{proof}
If such a solution exists, then it satisfies the Rankine-Hugoniot conditions. Thus $\sigma = \dfrac{f(u(0^-))-f(0)}{u_-}$, $z$ is continuous at $0$, and $-\sigma z'(x)=-kz(x)$ on $\R_-$, $z(x)=z(0)=1$ on $\R_+$ (as $z'(x)=0$ for all $x>0$), $u_0(x)=\lim_{y \rightarrow 0^+}u_0(y)=0$ for all $x>0$ (as, also, $(\sigma - f'(\uu(x))\uu'(x))=0$ for $x>0$) and $\dfrac{z_0(x)}{f'(u_0(x))-\sigma}=u_0'(x)$ on $\R^-$. \\

The conditions after are necessary as $f(\uu)'=\sigma \uu' +kq \uz$, and so by integrating we get $f(u_-)-\lim_{x\rightarrow -\infty}f(\uu(x))=\sigma (u_- - \lim_{x\rightarrow -\infty}\uu(x))+q\sigma$. This concludes the proof of the first point. \\

For the second point, as $f(\uu)-\sigma \uu$ is increasing on $\R_-$, and, for $x<0$ we have $\uu(x) >u_{-\infty}$, as otherwise there would exist $x_0$ such that $\uu(x_0)=u_{-\infty}$ and so $\sigma q =\int_{x_0}kq\uz(x)dx<\sigma q$. Hence, the solution to the ODE is globally defined and gives rise to a solution to the initial problem. This solution satisfies $\lim_{x\rightarrow -\infty}\uu(x)=u_{-\infty}$. \end{proof}
Note that, when the inverse function of $x \mapsto f(x)-\sigma x$ is known, one can obtain an explicit expression for $\uu$. \\

We will study the stability of such waves according to the sign of the propagation speed $\sigma$. The case $\sigma > 0$ will be studied first to obtain a stability result in a weighted space. The case $\sigma < 0$ will give rise to instability results in various spaces. The case of $\sigma = 0$ is not treated here. In fact, even the local in time existence is not contained in the framework of Theorem \ref{th1}. \\

We now can start to give some bounds on the wave that will prove useful later on. \\

Consider for the moment the case of speed $\sigma$ positive. We fix $k$, $f$ and $\uu(0^-)$ (such that $f'(\uu(0^-))>\dfrac{f(\uu(0^-))}{\uu(0^-)}>0$), and assume that $q \in \R$ is such that we can apply Proposition \ref{exstwave}. There is always such a $q$ as there is a sufficiently small $q_0 > 0$ such that the wave is defined for every $q \in [-q_0,q_0]$. \\

Let $\gamma > 0$ such that $f'>\sigma$ on $[\min(u_0,u_{-\infty})-\gamma,\max(u_{-\infty},u_0)+\gamma]=:I$.
We notice that for all $ x<0$ $$|\uu'(x)|\leq kqA\exp(k\sigma^{-1}x).$$
Furthermore, as $$\uu''(x)=\left( \dfrac{k}{\sigma}-f''(\uu(x))\dfrac{\uu'(x)}{(f'(\uu(x))-\sigma)} \right) \uu'(x),$$
and $$\uu^{(3)}(x)= \left( \dfrac{k}{\sigma}-\dfrac{f''(\uu(x))\uu'(x)}{(f'(\uu(x)-\sigma)} \right) ^2\uu'(x) - \dfrac{f^{(3)}(\uu(x))\uu'(x)^3 + f''(\uu(x))\uu''(x)\uu'(x)}{(f'(\uu(x))-\sigma)} + \dfrac{f''(\uu(x))^2\uu'(x)^2}{(f'(\uu(x))-\sigma)^2}.$$
Thus there exists some $\kappa=\kappa(q)$ real valued and defined on some neighborhood of $0$ such that $$|\uu'(x)| + |\uu''(x)| + |\uu^{(3)}(x)| \leq \kappa(q) \exp \left( \dfrac{kx}{\sigma} \right)$$ for all $q \in \R$ and $x \in \R^-$, and $\kappa$ being continuous with limit $0$ at $0$. \\
Similarly, we obviously have on $\R^-$ $$\uz^{(l)}(x)=\dfrac{k^l}{\sigma^l}\exp \left(\dfrac{kx}{\sigma} \right).$$ \\
\subsection{Perturbations}
Now, we are interested in damping estimates on solutions which are small perturbations to the original traveling waves. Until section 2.6, we assume $\sigma > 0$. \\

Thanks to the local wellposedness result in Theorem \ref{th1}, we obtain that there exists $\rho>0$ such that, given $(v_0,\zeta_0) \in H^{\,5/2}(\R \setminus \lbrace 0 \rbrace)$, supported away from the shock, with $L^{\infty}$ norm smaller than $\rho$. There exists $T \in (0,+\infty]$ such that there exists a unique maximal solution to the system with initial data $(\uu+v_0,\uz+\zeta_0)$. Remember that the shock localization $\psi$ is given by the Rankine-Hugoniot condition, $\psi(0)=0$, and for all $t\in [0,T)$ $\psi'(t)=\dfrac{f(u(t,\psi(t)^-)-f(u(t,\psi(t)^+)}{u(t,\psi(t)^-)-u(t,\psi(t)^+)}$. Furthermore,  $T<+\infty$ implies that we are at least in one of the following cases: a) the $L^{\infty}$ norm of $u(t,\cdot + \psi(t))-\uu$ or of $z(t,\cdot+\psi(t))-\uz$ becomes bigger than $\rho$ b) the space derivative of one of $u(t,\cdot+\psi(t))-\uu$ or of $z(t,\cdot+\psi(t))-\uz$ blows up. \\

We let $\zeta(t,x):=z(t,\psi(t)+x)-\uz(x)$ and $v(t,x):=u(t,x+\psi(t))-\uu(x)$ for all $(t,x) \in [0,T)\times \R^*$. We aim to show that finite time blow-up is prevented by choosing the initial perturbation sufficiently small, in the sense of some weighted Sobolev norm, and that the perturbation $(v,\zeta)$ goes to $0$ a $t$ goes to $+\infty$ in the weighted norm. \\

To do so, we will build an energy equivalent to the $H^2_{\eps}$ norm that will be non-increasing in $t$, by first introducing modifications of the norms on $L^2_{\eps}(\R^+)$ and $L^2_{\eps}(\R^-)$ of the form $E_{\pm}(w):=\int_{\R^{\pm}}w^2\varrho$ with the choices depending on the half-line considered ($\R_+$ or $\R_-$) and whether we consider terms in $\zeta$ or $v$. More precisely, we will consider them of the following forms $\varrho_{-,1}(x):=\exp \left(-\eps x-\int_0^x Ce^{-\eps |s|}ds \right)$ where $C$ is to be chosen later (taking $C$ big enough), $\varrho_{-,2}(x)=\exp \left( -\eps x + \int_0^x Ce^{-\eps |s|}ds \right)$ and $\varrho_+(x)=\exp(\eps x)$.\\

After that, the energy will be built as a sum of the terms of the form: $C_1(k,\pm)E_{\pm,1}(\partial^k_xv(t, \cdot))$ and $C_2(k,\pm)E_{\pm,2}(\partial^k_x\zeta(t, \cdot))$ where $k \in \lbrace 0,1,2 \rbrace$, and the $C_{\cdot}(k,\pm)$ are constants independent of $(v_0,\zeta_0)$. \\

Here are the main ingredients in our proof of the existence of such an energy function, and also how fast it goes to $0$: \\

The choices of the constants $C_{j}(k,\pm)$ will be made to control the terms that, when taking the time derivatives of the $E_{\pm}$ of $v$, $\zeta$ or one of their spatial derivatives up to order $2$, is non-negative. Terms of the form $\partial^l_x\zeta(t,0^-)$ that will appear due to the outgoing characteristic, will be bounded through boundary terms appearing in integration by parts done on integrals appearing in the derivatives of the $E_+(\partial^s_x\zeta(t, \cdot))$ for $s \leq l$, and terms that have the form $|\psi'(t)-\sigma|$ will be controlled through the boundary terms of $E_{\pm}(v(t, \cdot))$. \\

Finally, the choice of a perturbations of the norms of $L^2_{\eps}(\R^+)$ and $L^2_{\eps}(\R^-)$ by terms of the form $\exp\left( \pm \int_0^xe^{-\eps |s|}ds )\right)$ is made to bound some terms that are not small, but are of integral form with a quadratic term in $v$, $\zeta$ and some of their spatial derivatives times some derivative of the underlying shock profile. \\

Now, we can try to obtain bounds on the size of the perturbations $(v,\zeta)$, first by writing done the equations they satisfy, as well as the equations their derivatives satisfy. \\

From now on, we assume the initial data to be a smooth and compactly supported function in $\R^*$. With $T$, its (possibly infinite) time of existence, for every time $t \in [0,T)$, the solution is in $H^s(\R^*)$ for every $s \in \R_+$, and, as long as the $L^{\infty}$ norm of $v$ is small enough and the solution remains bounded in $W^{1,\infty}(\R^*)$, we have on $\R^-$ 
$$v_t(t, \cdot)=(\psi'(t)-f'(\uu+v(t, \cdot)))v_x(t, \cdot)+(\psi'(t)-\sigma +f'(\uu)-f'(\uu+v(t, \cdot)))\uu'+kq\zeta(t, \cdot),$$
and $$\zeta_t=(\psi'(t)-\sigma)\uz' -k\zeta(t, \cdot) + \psi'(t)\zeta_x(t, \cdot).$$
Whereas \begin{align*}
v_{tx}(t,.)=&(\psi'-f'(\uu+v(t, \cdot)))v_{xx}(t, \cdot)-f''(\uu+v(t, \cdot))(\uu'+v_x(t, \cdot))v_x(t, \cdot) \\ 
&+ (\psi'-\sigma+f'(\uu)-f'(\uu+v(t, \cdot)))\uu'' +(f''(\uu)-f''(\uu+v(t, \cdot)))\uu'^2 \\ &-f''(\uu+v(t, \cdot))\uu'v_x(t, \cdot) +kq\zeta_x(t, \cdot), 
\end{align*}
\begin{align*}
v_{txx}(t, \cdot)=&(\psi'-f'(\uu+v(t, \cdot)))v_{xxx}(t, \cdot)-2f''(\uu+v(t, \cdot))(\uu'+v_x(t, \cdot))v_{xx}(t, \cdot) \\ 
&-f^{(3)}(\uu+v(t, \cdot))(\uu'+v_x(t, \cdot))^2v_x(t, \cdot) - f''(\uu + v(t, \cdot))(\uu'' + v_{xx}(t, \cdot))v_{xx}(t, \cdot) \\ & +(\psi'-\sigma+f'(\uu)-f'(\uu+v(t, \cdot)))\uu^{(3)} + 2(f''(\uu)-f''(\uu+v(t, \cdot)))\uu'\uu''  \\ & - 2f''(\uu+v(t, \cdot))\uu''v_x(t, \cdot) + (f''(\uu)-f''(\uu+v(t, \cdot)))\uu'\uu'' \\ & + (f^{(3)}(\uu)-f^{(3)}(\uu+v(t, \cdot)))\uu'^3 - f^{(3)}(\uu+v(t, \cdot))(2\uu'+v_x(t, \cdot))\uu'v_x(t, \cdot) \\ & - f''(\uu+v(t, \cdot))\uu'v_{xx}(t, \cdot)+kq\zeta_{xx}(t, \cdot),
\end{align*}
\begin{align*}
\zeta_{tx}(t, \cdot)=&(\psi'(t)-\sigma)\uz''-k\zeta_x(t, \cdot)+\psi'(t)\zeta_{xx}(t, \cdot), \\
\zeta_{txx}(t, \cdot)=&(\psi'(t)-\sigma)\uz^{(3)}-k\zeta{xx}(t, \cdot)+\psi'(t)\zeta_{xxx}(t, \cdot).
\end{align*}
Similarly, on $\R^+$ $$v_t(t, \cdot)=(\psi'(t)-f'(v(t, \cdot)))v_x(t, \cdot),$$ and $$\zeta_t(t, \cdot)=\psi'(t)\zeta_{tx}(t, \cdot).$$
Furthermore \begin{align*}
v_{tx}(t, \cdot)=&(\psi'(t)-f'(v(t, \cdot)))v_{xx}(t, \cdot)-f''(v(t, \cdot))(v_x(t, \cdot))v_x(t, \cdot), 
\end{align*}
\begin{align*}
v_{txx}(t, \cdot)=&(\psi'(t)-f'(v(t, \cdot)))v_{xxx}(t, \cdot)-2f''(v(t, \cdot))(v_x(t, \cdot))v_{xx}(t, \cdot) \\ &-f''(v(t, \cdot))(v_{xx}(t, \cdot))v_x(t, \cdot)-f^{(3)}(v(t, \cdot))(v_x(t, \cdot))^2v_x(t, \cdot), 
\end{align*}
\begin{align*}
\zeta_{tx}(t, \cdot)=&\psi'\zeta_{xx}(t, \cdot), 
\end{align*}
\begin{align*}
\zeta_{txx}(t, \cdot)=&\psi'(t)\zeta_{xxx}(t, \cdot).
\end{align*}
\subsection{Boundary terms}
In the energy estimates, we will need to control boundary terms that can not be neglected due to the characteristics outgoing from the shock, specifically those involving $\zeta$ or one of its derivatives at $0^-$, by using other boundary terms going into the shock, in particular only those involving $\zeta$ or one of its derivatives at $0^+$. We also need to control the derivatives of the phase by using the boundary terms involving $v$ or one of its derivatives both at $0^+$ and $0^-$, where we note both of them are linked to characteristics going into the shock. \\

To obtain such bounds, we will use the Rankine-Hugoniot condition, and the fact that for smooth enough solutions we can differentiate it with respect to time, and, using the equation, replace the time derivatives of $v$ and $\zeta$ with terms involving only spatial derivatives.
\begin{lemma}
If $\eta >0$ is small enough, then there exists a constant $\tilde{C_f}$ such that, for $T>0$ and every solution $(v,\zeta)$ (either for a fixed $q$ or for uniformly in $q$ for $q$ small) initially a smooth function compactly supported away from $0$, such that $(v(t, \cdot),\zeta(t, \cdot))$ is defined on $[0,T]$ with its $L^{\infty}$ norm smaller than $\eta$, then we have, on $[0,T]$ $|\psi''(t)| \leq \tilde{C}_f \|(v,\zeta)\|_{W^{1,\infty}}$, $| \psi ' (t) - \sigma |^2 \leq \tilde{C}_f (|v(t,0^-)|^2+|v(t,0^+)|^2)$.  \\

And such that we also have, on $[0,T]$
	$$\zeta(t,0^-)=\zeta(t,0^+),$$
 	$$\zeta_x(t,0^-)^2 \leq \tilde{C}_f(\zeta_x(t,0^+))^2 + \tilde{C}_fk\zeta(t,0^-)^2+\tilde{C}_f(v(t,0^+)^2+v(t,0^-)^2),$$ and 
 	\begin{align*}
 		(\zeta_{xx}(t,0^-))^2 \leq& \tilde{C}_f(\zeta_{xx}(t,0^+))^2 + \tilde{C}_f(\zeta_x(t,0^+))^2 + \tilde{C}_f(v(t,0^-)^2+v(t,0^+)^2)\\
 		& + \tilde{C}_f\zeta(t,0^+)^2 + \tilde{C}_f(v_x(t,0^-)^2 + v_x(t,0^+)^2)
 	\end{align*}
\end{lemma}
\begin{proof}We impose $\eta \leq \dfrac{u_0}{4}$. For the bounds on the values of $\zeta$ and its spatial derivatives at $0^-$, we will use that $$\zeta(t,0^+)=\zeta(t,0^-)$$ for every $t$, and so that $$\partial^l_t\zeta(t,0^+)=\partial^l_t\zeta(t,0^-)$$ for every $t$ and $l$. \\

$$\zeta_t(t,0^+)=\psi'(t)\zeta_x(t,0^+),$$
$$\zeta_t(t,0^-)=\psi'(t)\zeta_x(t,0^-)-k\zeta(t,0^-)+(\psi'-\sigma)\dfrac{k}{\sigma},$$
$$\zeta_{tt}(t,0^+)=\psi''(t)\zeta_x(t,0^+)+\psi'(t)\zeta_{tx}(t,0^+),$$
$$\zeta_{tt}(t,0^-)=\psi'(t)\zeta_{tx}(t,0^-)+\psi''(t)\zeta_x(t,0^-)-k\zeta_t(t,0^-)+\psi''(t)\dfrac{k}{\sigma}.$$
We also have that $$\psi''(t)=-\dfrac{\psi'(t)(v_t(t,0^-)-v_t(t,0^+))}{u_0+v(t,0^-)-v(t,0^+)}+\dfrac{f'(u_0+v(t,0^-))v_t(t,0^-)-f'(v(t,0^+))v_t(t,0^+)}{u_0+v(t,0^-)-v(t,0^+)}.$$
And so \begin{align*}
|\psi''(t)|\leq& \dfrac{1}{u_0 - 2\eta} ( \nu ( k|q|\zeta(t,0^-)| + C_f^2k|q|(|v(t,0^-)|+|v(t,0^+)|) + \nu (|v_x(t,0^-)| + |v_x(t,0^+)|)) \\ &+ C_f ( \nu ( k|q|\zeta(t,0^-)| + C_f^2k|q|(|v(t,0^-)|+|v(t,0^+)|) + \nu (|v_x(t,0^-)| + |v_x(t,0^+)|))).
\end{align*}
We obtain that $$\zeta_x(t,0^-)=\zeta_x(t,0^+)+\dfrac{k\zeta(t,0^-)}{\psi'(t)}+\dfrac{(\psi'(t)-\sigma)k}{\sigma\psi'(t)}.$$
Using that, for all $(a,b,c)\in \mathbb{R}$, $(a+b+c)^2\leq 3(a^2+b^2+c^2)$, we have \begin{align*}
\zeta_x(t,0^-)^2 \leq& 3(\zeta_x(t,0^+))^2 + 3\left( \dfrac{k\zeta(t,0^-)}{\mu}\right) ^2+3 \dfrac{\tilde{C}_f^2(v(t,0^+)^2+v(t,0^-)^2)k^2}{\sigma^2\mu^2}.
\end{align*}
We have, for second order derivative that \begin{align*}
\psi'(t)\zeta_{xx}(t,0^-)=&+k\zeta_x(t,0^-)-(\psi(t)-\sigma)\dfrac{k^2}{\sigma ^2} + \psi'(t)\zeta_{xx}(t,0^+) + k\zeta_x(t,0^+) + \dfrac{k^2\zeta(t,0^-)}{\psi'(t)} \\ & + \dfrac{k(\psi'(t)-\sigma)}{\psi'(t)}\dfrac{k}{\sigma} - \dfrac{\psi''(t)\zeta(t,0^-)}{\psi'(t)^2} + k\dfrac{\psi''(t)}{\psi'(t)^2},
\end{align*}
Thus
\begin{align*}
(\zeta_{xx}(t,0^-))^2 \leq& \tilde{C}_f(\zeta_x(t,0^+))^2 + \tilde{C}_fk\zeta(t,0^-)^2+\tilde{C}_f(v(t,0^+)^2+v(t,0^-)^2).
\end{align*}
\end{proof}	
As $\psi'(t) \in [\sigma-C_f\eta,\sigma+C_f\eta]$ and $f'(\uu+v)\in [-C_f\eta,C_f\eta]\cup f'([u_{-\infty}-C_f\eta,u_0+C_f\eta])$, we have, by shrinking $\eta$ if necessary, that there exists $\mu > 0$ such that $\psi'(t)>\mu$ for every $x>0$, $\psi'(t)-f'(\uu(x)+v(x))>0$, and for every $x<0$, $\psi'(t)-f'(\uu(x)+v(x))<-\mu$, and $k \geq \mu$, and there exists $\nu>0$ such that $\psi'(t)<\nu$.  \\
\subsection{Energy estimates}
We will now obtain estimates of the time derivatives of the $E_{\pm}$ of $v$ and $\zeta$ and their spatial derivatives. Below, $\delta_1$ and $\delta_2$ will denote some positive constants to be determined later. At the end, we will choose two values, one to control terms in $v$ or $v_x$ or $v_{xx}$ and the other to control terms in $\zeta$ or $\zeta_x$ or $\zeta_{xx}$. \\

Furthermore, we will derive estimates for $\uu(0^-)$, $f$, $k$ fixed while $q$ is a parameter. Depending on $q$, we will obtain either a high-frequency damping estimate, or, a direct nonlinear stability result in the weak heat release limit. \\

The idea, which we will soon check directly at the linear $L^2$ level, is that by choosing the constant $C$ big enough and $q$ be smaller (in absolute value) than some $q_0 > 0$, there exists three positive constants $C_2$, $C_3$ and $C_4$ as well as some $\nu > 0$ (independent of $v$ and $\zeta$) such that if $\|v(t,\cdot)\|_{W^{1,\infty}}$ and $\|\zeta(t,\cdot)\|_{W^{1,\infty}}$ are smaller than $\eta$ on $[0,T]$, then, on $[0,T]$ \begin{align*}\dfrac{d}{dt}(E_{-,1}(v)+C_2E_{-,2}(\zeta)&+C_3E_+(v)+C_4E_+(\zeta))(t,\cdot) \leq \\ &-\nu (E_{-,1}(v(t,\cdot))+C_2E_{-,2}(\zeta(t,\cdot))+C_3E_+(v(t,\cdot))+C_4E_+(\zeta(t,\cdot))).
\end{align*} 
From here, the idea is to obtain similar bounds for the higher order terms, while taking into account that there may be some loss of lower order derivatives when controlling higher order terms. \\

The linear problem is given by \begin{align*}
v_t=&  (\psi' - \sigma)\uu' + kq \zeta- ((f'(\uu) - \sigma)v)_x  \quad \text{on} \, \, \R^-, \\
\zeta_t=&\sigma \zeta_x - k \zeta + (\psi' - \sigma)\uz' \quad  \text{on} \, \, \R^-, \\
\zeta_t=&\sigma \zeta_x \quad  \text{on} \, \, \R^+, \\
v_t=&(\sigma - f'(0))v_x \quad  \text{on} \, \, \R^+, \\
\psi'=&\dfrac{(f'(u_0)-\sigma)v(\cdot,0^-)}{u_0} + \dfrac{(\sigma - f'(0))v(\cdot,0^+)}{u_0}, \; \; \; \; \zeta(\cdot,0^+)=\zeta(\cdot,0^-).
\end{align*}
In fact, we can obtain, by fixing several positive parameters $C$, $\tilde{\delta}$, $\eta$ and $\eps:=\dfrac{k}{\nu}$ that \begin{align*}
&\dfrac{d}{dt}\left( E_{-,1}(v) + E_+(v) + E_+(\zeta) +  E_{-,2}(\zeta) \right) \leq \\ &v(t,0^-)^2 \left( \dfrac{\sigma - f'(u_0)}{2} + \dfrac{C_f \tilde{\delta}}{2} \int_{\R_-}|\uu'|e^{-\eps \cdot - C\int_0^{\cdot}} +  C_f\dfrac{\eta}{2}\int_{\R_-}\uz'e^{-\eps + \int_0^{\cdot}Ce^{\eps s}ds} \right) \\ & + v(t,0^+)^2 \left(  \dfrac{f'(0)-\sigma}{2} + \dfrac{C_f \tilde{\delta}}{2} \int_{\R_-}|\uu'|e^{-\eps \cdot - C\int_0^{\cdot}} +  C_f\dfrac{\eta}{2}\int_{\R_-}\uz'e^{-\eps + \int_0^{\cdot}Ce^{\eps s}ds} \right) \\ &+ \int_{\R_-}\zeta^2 \left( \dfrac{k|q|}{2\delta} e^{2C\eps^{-1}} +  \dfrac{\eps \sigma - 2k}{2} + \dfrac{(\eta^{-1} - C)}{2}e^{\eps \cdot} \right) e^{-\eps \cdot + C \int_0^{\cdot} e^{\eps s}ds} 
\end{align*}
\begin{align*}
 &+ \int_{\R_-}v^2 \left( -\dfrac{C_f\eps}{2} + \dfrac{k|q|\delta}{2} + \dfrac{(k |q| C_f \tilde{\delta}^{-1} + C_fk|q| - C)e^{\eps \cdot}}{2}\right) e^{-\eps \cdot - C\int_0^{\cdot}e^{\eps s}ds}  \\ &+ \int_{\R_+}\zeta^2 \left( - \dfrac{\eps }{\sigma} \right) e^{\eps \cdot} + \int_{\R_+}v^2 \left( - \dfrac{(f'(0)-\sigma)}{2} \right) e^{\eps \cdot}
\end{align*}
We thus have it is enough to have
\begin{itemize}
	\item $\dfrac{(f'(0)-\sigma)}{2} + \dfrac{C_f\sigma \tilde{\delta}e^{C\eps^{-1}}}{k} + \dfrac{ C_f \eta \sigma}{k} \leq 0$ (constraint associated to $v(t,0^+)^2$), 
	\item $\dfrac{k|q| e^{2C\eps^{-1}}}{2\delta} +  \dfrac{\eps \sigma - 2k}{2} + \dfrac{(\eta^{-1} - C)}{2} < 0$ (constraint associated to $\zeta^2$),
	\item $\dfrac{k|q|\delta}{2} - \dfrac{C_f \eps}{2} < 0$ (term in $v^2$ on $\R^-$, first part),
	\item $f'(u_0) - \sigma \geq 2 C_f^2e^{C\eps^{-1}} \tilde{\delta} + C_f\eta$ (term in $v(t,0^-)^2$),
	\item $k|q|\tilde{\delta}^{-1} + C_f k |q| - C \leq 0$ (term in $v^2$ on $\R^-$, second part).
\end{itemize}

We now choose the constants as follows: $\eta \leq \dfrac{(\sigma - f'(0))k}{4C_f\sigma}$ (as it needs to satisfy other conditions listed before, all independent on $q$ crucially), $C = \eta^{-1}$, $\tilde{\delta} = \dfrac{\sigma - f'(0)}{C_f \sigma e^{2C\eps^{-1}}}k$, giving us the result we were aiming for as long as $|q|$ is small enough. \\

The computations in this section are straightforward but written here for the sake of completeness. To state our goal, it is to isolate boundary terms coming from integration by parts in order to reduce the order of the derivatives appearing in the integrals; thus to only have derivatives of order at most $k$ in the bound of the time derivative of $E_{\pm}(\partial^k_xw)$ where $w$ is $v$ or $\zeta$, as well as isolating terms that cannot be controlled directly by some form of dissipation (in particular, terms related to the phase or the reaction term when computing the time derivative of $E_-(v(t, \cdot))$), and, if possible, have the uncontrolled terms multiplied either by some small constant (that will depend on the energy studied being the one of $v$ or one of its derivatives, or of $\zeta$ or of one of its derivatives) or by some exponentially localized terms in space (inside of the integral). \\

We have for any $\eps > 0$, $C > 0$, $\delta_1 > 0$ and as long as $\sup_{0\leq t\leq T} \|(v(t, \cdot),\zeta(t, \cdot))\|_{W^{1,\infty}(\R^*)} < \eta$, the following inequalities hold for all $t$
\begin{align*}
\dfrac{1}{2}\dfrac{dE_{-,1}(v)}{dt}(t)\leq& \dfrac{k|q|}{2\delta_1}E_{-,2}(\zeta(t,\cdot)) \exp \left( \dfrac{2C}{\eps} \right)+\dfrac{k|q|\delta_1}{2}E_{-,1}(v(t,\cdot)) - \mu \dfrac{v(t,0^-)^2}{2} \\ & - \int_{\R^-}\mu ( \eps + C e^{\eps \cdot} )\dfrac{v(t,\cdot)^2}{2}\varrho_{-,1}dx +\int_{\R^-}C_f(|\uu'|+\eta)\dfrac{v(t,\cdot)^2}{2}\varrho_{-,1}dx \\ & + \dfrac{1}{2} | \psi'(t) - \sigma |^2 E_{-,1}(\sqrt{\uu'}) + \dfrac{1}{2}E_{-,1}(\sqrt{\uu'}v(t,\cdot)) ,
\end{align*}
\begin{align*}
\dfrac{1}{2}\dfrac{dE_{-,1}(v_x)}{dt}(t)\leq& \int_{\R_-}(v_x(t,\cdot))^2\left( -\dfrac{\mu \eps + C \mu e^{\eps \cdot}}{2} + \dfrac{C_f}{2}\eta + \dfrac{C_f}{2}|\uu'| + \dfrac{k|q|\delta_1}{2} + \dfrac{|\uu'|}{2} \right. \\ & \left. + \dfrac{C_f(\uu'^2 + |\uu''|)}{2} + |\uu'|C_f \right) \varrho_{-,1} + \int_{\R_-}k|q|\delta_1v(t,\cdot)^2\dfrac{\varrho_{-,1}}{2} \\ & + \dfrac{k|q|e^{2C\eps^{-1}}E_{-,2}(\zeta_x)}{2\delta_1} + \dfrac{|\psi'-\sigma|^2E_{-,1}(\sqrt{|\uu'|})}{2}
\end{align*}
\begin{align*}
\dfrac{1}{2}\dfrac{dE_{-,1}(v_{xx}(t,\cdot))}{dt}(t) \leq& \int_{\R_-}(v_{xx}(t,\cdot))^2\left( -\dfrac{\mu \eps + \mu C e^{\eps \cdot}}{2} + |\uu'| \left( \dfrac{(5+2\eta)C_f + C_f|\uu'| + 3|\uu''| + \uu'^2}{2} \right) \right. \\ & \left.+ \dfrac{5C_f}{2}|\uu''|+ \dfrac{C_f + 1}{2} |\uu^{(3)}| + \dfrac{7C_f\eta + \delta_1k|q| + C_f\eta^2}{2} \right)\varrho_{-,1} \\ & + \int_{\R_-}(v_x(t, \cdot))^2 \left( \dfrac{C_f ( 3|\uu''| + 3 \eta |\uu'| + \eta^2 + \uu'^2 )}{2} \right) \varrho_{-,1} \\ & + \int_{\R_-} v(t, \cdot)^2 \left( C_f \dfrac{|\uu^{(3)}| + 3|\uu'\uu''| + |\uu'|^3}{2} \right) \varrho_{-,1} \\ & + |\psi'-\sigma|^2 \dfrac{E_{-,1}(\sqrt{|\uu^{(3)}|})}{2} + \dfrac{k|q|e^{2C\eps^{-1}}}{2\delta_1} E_{-,2}((\zeta_{xx}(t, \cdot))) . 
\end{align*}
\begin{align*}
\dfrac{1}{2}\dfrac{dE_{-,2}(\zeta)}{dt}(t) \leq&  \dfrac{1}{2} | \psi ' (t) - \sigma |^2E_-(\sqrt{\uz'}) + \dfrac{1}{2}E_{-,2}(\sqrt{\uz'} \zeta^2) + \nu \dfrac{\zeta(0^-)^2}{2} \\ & +  \int_{\R^-}\nu \eps \dfrac{\zeta^2}{2}\varrho_{-,2}dx- \int_{\R_-} C \mu e^{\eps \cdot}\dfrac{\zeta^2}{2}\varrho_{-,2} - k E_{-,2}(\zeta),
\end{align*}
\begin{align*}
\dfrac{1}{2}\dfrac{dE_{-,2}(\zeta_x)}{dt}(t)\leq& \nu\dfrac{(\zeta_x)(0^-)^2}{2}+\nu \eps \int_{\R^-}\dfrac{(\zeta_x)^2}{2}\varrho_{-,2} -\int_{\R^-}C\nu e^{\eps \cdot}\dfrac{(\zeta_x)^2}{2}\varrho_{-,2} - kE_{-,2}(\zeta_x) \\ &+ \dfrac{| \psi ' (t) - \sigma |^2 E_{-,2}(\sqrt{|\uz''|})}{2} + \dfrac{E_{-,2}(\sqrt{|\uz''|}\zeta_x)}{2}, 
\end{align*}
\begin{align*}
\dfrac{1}{2}\dfrac{dE_{-,2}(\zeta_{xx})}{dt}(t)\leq& \dfrac{\nu\zeta_x(0^-)^2}{2} + \int_{\R^-}\dfrac{(\zeta_x)^2}{2}(\nu \eps - C \mu e^{\eps \cdot})\varrho_{-,2} - k E_{-,2}(\zeta_{xx})
\\& + \dfrac{| \psi ' (t) - \sigma |^2 E_{-,2}(\sqrt{|\uz^{(3)}|})}{2} + \dfrac{E_{-,2}(\sqrt{|\uz^{(3)}|}\zeta_{xx})}{2}, 
\end{align*}
\begin{align*}
\dfrac{1}{2}\dfrac{dE_+(\zeta)}{dt}(t) \leq - \mu \dfrac{\zeta (t,0^+)^2}{2} - \mu \int_{\R^+} \eps \dfrac{\zeta^2}{2}\varrho_+, 
\end{align*}
\begin{align*}
\dfrac{1}{2}\dfrac{dE_+(\zeta_x)}{dt}(t)\leq-\mu\dfrac{\zeta_x(t,0^+)^2}{2}-\mu\eps\int_{\R^+}\dfrac{(\zeta_x)^2}{2}\varrho_+, 
\end{align*}
\begin{align*}
\dfrac{1}{2}\dfrac{dE_+(\zeta_{xx})}{dt}(t)\leq-\mu\dfrac{\zeta_{xx}(t,0^+)^2}{2}-\mu\eps\int_{\R^+}\dfrac{(\zeta_{xx})^2}{2}\varrho_+,
\end{align*}
\begin{align*}
\dfrac{1}{2}\dfrac{dE_+(v(t, \cdot))}{dt}&\leq -\dfrac{\mu v(t,0^+)^2}{2} - \int_{\R_+} \mu \dfrac{v(t, \cdot)^2}{2}(\eps + Ce^{\eps \cdot}) \varrho_+ + \dfrac{C_f\eta}{2} \int_{\R_+}v(t, \cdot)^2\varrho_+,
\end{align*}
\begin{align*}
\dfrac{1}{2}\dfrac{dE_+(v_x)}{dt}(t)\leq&-\mu\dfrac{(v_x(t,0^+))^2}{2}-\int_{\R^+}\mu\eps(v_x(t, \cdot))^2\dfrac{\varrho_+}{2}+\int_{\R^+}C_f\dfrac{\eta(v_x(t, \cdot))^2}{2}\varrho_+,
\end{align*}
\begin{align*}
\dfrac{1}{2}\dfrac{dE_+(v_{xx})}{dt}(t) \leq& \dfrac{5C_f \eta - 2 \mu \eps}{4} \int_{\R^+}(v_{xx}(t, \cdot))^2 + \int_{\R^+}\dfrac{-C \mu e^{\eps \cdot}(v_{xx}(t, \cdot))^2}{2}\varrho_+ \\ & - \dfrac{\mu v_{xx}(t,0^+)^2}{2} + \int_{\R^+} \dfrac{C_f}{2} \eta^2((v_x(t, \cdot))^2+(v_{xx}(t, \cdot))^2)\rho_+, 
\end{align*}
Here, $\mu$ is the same parameter from the end of the proof of Proposition \ref{exstwave}. Thus, we now assume $\eta > 0$ satisfies $\eta \leq \min\left( \dfrac{u_0}{4}, 1, \dfrac{\mu\eps}{32} \right)$. We will impose more conditions related to the well-posedness result Theorem \ref{th1} on $\eta$ later on.
When it comes to $\delta_1$, we choose it such that $k|q|\delta_1 \leq \dfrac{\mu \eps}{8}$ (more precisely, we take $\delta_1 = 1$ if $q=0$ and $\delta_1=\dfrac{\mu\eps}{8k|q|}$ otherwise). The following choice is made for $C$ \begin{align*}
C=2\mu^{-1}\max &\left(\kappa(q)\dfrac{11C_f + 2\eta C_f + \kappa(q) (C_f + 3 + \kappa(q))}{2} + \dfrac{1 + \eta(7+\eta) + \delta_1k|q|}{2} , \right. \\
& \; \;\, \left.  \kappa(q) \dfrac{4C_f+C_f\kappa(q)+1}{2}, \kappa(q) \dfrac{C_f + 1}{2} , \dfrac{k+k^2+k^3}{\nu} \right)
\end{align*} and define $\omega := \dfrac{\mu \eps}{4}.$ \\

When it comes to choosing $q_0$, we first fix the other constants by first taking some $q_1 > 0$ for which the $\uu$ and $\uz$ are defined for every $q \in [-q_1,q_1]$ and for every such $q$ with $|q|\leq q_1$ $\kappa(q) \leq 1$, and then replacing in the definition of constants the ones obtained with $\kappa(q)$ replaced by $1$ and $q$ by $q_1$. Finally, $q_0$ will be chosen later on. \\

There exists a constant $\tilde{C}$ such that, as long as the $W^{1,\infty}(\R^*)$ norm of the perturbation remains smaller than $\eta$, we have \begin{align*}
(\zeta(t,0^-))^2 \leq \tilde{C}(\zeta(t,0^+))^2
\end{align*}
\begin{align*}
(\zeta_x(t,0^-))^2 \leq \tilde{C}(\zeta_x(t,0^+)^2 + |v(t,0^+)|^2 + |v(t,0^-)|^2 + \zeta(t,0^+)^2)
\end{align*}
\begin{align*}
(\zeta_{xx}(t,0^-))^2 \leq \tilde{C}(v(t,0^+)^2 + v(t,0^-)^2 + v_x(t,0^+)^2 + v_x(t,0^-)^2 + \zeta(t,0^+)^2 + \zeta_x(t,0^+)^2)
\end{align*}
Thus, we can now adjust the value of $\tilde{C}$ in a way which only depends on $q_1$ (while $\tilde{C}$ depends on $k$, $f$ and $\uu(0^-)$, we fixed these quantities independently of $q$ and so we are not worried by that dependence) such that, for every such $q \in [-q_0,q_0]$ and solution $(v,\zeta)$ defined on $[0,T)$ which remains strictly smaller than $\eta$ on $[0,T']$ we have that, on $[0,T']$, if $q \in [-q_0,q_0]$ \begin{align*}
\dfrac{dE_{-,1}(v)(t)}{dt} \leq& -\omega E_{-,1}(v(t, \cdot)) + \tilde{C}|q|E_{-,2}(\zeta(t, \cdot)) - \omega v(t,0^-)^2 + \tilde{C}(v(t,0^+))^2\sqrt{\kappa(q)}, \\
\dfrac{dE_{-,1}(v_x)(t)}{dt} \leq&  - \omega v_x(t,0^-)^2 + \tilde{C}E_{-,1}(v(t, \cdot)) + \tilde{C}((v(t,0^-))^2 + (v(t,0^+))^2)\sqrt{\kappa(q)} - \omega E_{-,1}(v_x(t, \cdot)) \\ &  + \tilde{C}|q|E_{-,2}(\zeta_x(t, \cdot)), \\
\dfrac{dE_{-,1}(v_{xx})(t)}{dt} \leq&  - \omega v_{xx}(t,0^-)^2 + \tilde{C}(E_{-,1}(v(t, \cdot)) + E_{-,1}(v_x(t, \cdot))) + \tilde{C}((v(t,0^-))^2 + (v(t,0^+))^2)\sqrt{\kappa(q)}-\\& - \omega E_{-,2}(v_{xx}(t, \cdot)) + \tilde{C}|q|E_{-,2}(\zeta_{xx}(t, \cdot)), 
\end{align*}
\begin{align*}
\dfrac{dE_+(v)(t)}{dt} \leq& - \omega\dfrac{v(t,0^+)^2}{2} - \omega E_+(v(t, \cdot))
\end{align*}
\begin{align*}
\dfrac{dE_+(v_x)(t)}{dt} \leq& -\omega E_+(v_x(t, \cdot)) - \dfrac{\omega(v(t,0^+))^2}{2}, 
\end{align*}
\begin{align*}
\dfrac{dE_+(v_{xx})(t)}{dt} \leq& -\omega E_+(v_{xx}(t, \cdot)), 
\end{align*}
\begin{align*}
\dfrac{dE_{-,2}(\zeta)(t)}{dt} \leq& -\omega E_{-,2}(\zeta(t, \cdot)) + \tilde{C}(v(t,0^+)^2 + v(t,0^-)^2 + \zeta(t,0^+)^2), 
\end{align*}
\begin{align*}
\dfrac{dE_{-,2}(\zeta_x(t, \cdot))}{dt} \leq& -\omega E_{-,2}(\zeta_x(t, \cdot)) + \tilde{C}(v(t,0^+)^2+v(t,0^-)^2) + \zeta_x(t,0^+)^2 + \zeta(t,0^+)^2), 
\end{align*}
\begin{align*}
\dfrac{dE_{-,2}(\zeta_{xx})(t)}{dt} \leq& -\omega E_{-,2}(\zeta_{xx}(t, \cdot)) + \\&
 + \tilde{C}(v(t,0^+)^2+v(t,0^-)^2 + \zeta_{xx}(t,0^+)^2 + \zeta_x(t,0^+)^2 + \zeta(t,0^+)^2 + v_x(t,0^+)^2 + v_x(t,0^-)^2), 
\end{align*}
\begin{align*}
\dfrac{dE_+(\zeta)}{dt}(t) \leq& -\omega (E_+(\zeta(t, \cdot)) + \zeta(t,0^+)^2), 
\end{align*}
\begin{align*}
\dfrac{dE_+(\zeta_x)}{dt}(t) \leq& -\omega (E_+(\zeta_x(t, \cdot)) + \zeta_x(t,0^+)^2), 
\end{align*}
\begin{align*}
\dfrac{dE_+(\zeta_{xx})}{dt}(t) \leq& -\omega (E_+(\zeta_{xx}(t, \cdot)) + \zeta_{xx}(t,0^+)^2).
\end{align*}
\subsection{Nonlinear stability and high-frequency damping estimates}
We will now be able to conclude the argument in both cases. For stability, we will choose the constants in order to be able to obtain the decay in time of the energy and for the high-frequency damping estimates we choose the constants to control the rest of the equations by lower order terms. \\

We will first focus on the small $q$ behavior. \\
\begin{proposition}\label{prop:constants} There exists a $12$-tuple of positive constants $(C_{0,-},C_{1,-},C_{2,-},C_{0,-}',...)$, a positive constant $q_0$ small enough, $C>0$ big enough, $\vartheta >0$ such that, for every $q \in [-q_0,q_0]$ and every perturbation of the wave that initially satisfies the conditions of \ref{th1}, namely $(v_0,\zeta_0)$ supported away from zero and smooth, there exists a unique maximal solution $(v,\zeta)$ to the problem associated with $q$, with initial data $(\uu + v_0,\uz + \zeta_0)$ defined on some time interval $[0,T)$ with $T \in (0,+\infty]$. Then, for every $T' \leq T$ such that $\|(v,\zeta)\|_{W^{1,\infty}} < \eta$ for every $t \in [0,T')$, we have, on $[0,T')$
\begin{align*}
&\dfrac{1}{2}\sum_{k=0}^2(C_{k,-}\dfrac{dE_{-,1}(
	\d_x^k v)}{dt}(t)+C_{k,-}'\dfrac{dE_{-,1}(\d_x^k\zeta)}{dt}(t) + C_{k,+}\dfrac{dE_+(\d_x^kv)}{dt}(t) + C_{k,+}'\dfrac{dE_+(\d_x^k\zeta))}{dt}(t) \\ &\leq -\dfrac{\vartheta}{2}\sum_{k=0}^2(C_{k,-}E_{-,1}(\d_x^kv(t, \cdot)+C_{k,-}'E_{-,1}(\d_x^k\zeta(t, \cdot)) + C_{k,+}E_+(\d_x^kv(t, \cdot)) + C_{k,+}'E_+(\d_x^k\zeta(t, \cdot)))
)\end{align*} 
Thus, for some constant $M>0$ independent on $q$, $v_0$, $\zeta_0$ and $T$, for every $t \in [0,T)$, if $\|(v(s),\zeta(s)\|_{W^{1,\infty}} < \eta$ for every $s \in [0,t]$, then $$\sum_{k=0}^2(\|\partial^k_xv(t, \cdot)\|_{L_{\eps}^2(\mathbb{R}^*)}^2+\|\partial^k_x\zeta(t, \cdot)\|_{L_{\eps}^2(\R^*)}) \leq Me^{-(t-s)\vartheta}\sum_{k=0}^2(\|\partial^k_xv(s, \cdot)\|_{L_{\eps}^2(\mathbb{R}^*)}^2+\|\partial^k_x\zeta(s, \cdot)\|_{L_{\eps}^2(\R^*)})$$
\end{proposition}
\begin{proof}
To obtain the desired inequality, we just need to ensure that we can force certain coefficients to non-positive if $q_0$ is chosen small enough and $C > 0$ big enough. We will do this by using the following bounds \\

We may choose $C_{0,-}$, $C_{1,-}$, $C_{2,-}$, ... such that they satisfy the following inequalities \\ \begin{itemize}
\item(to have a nonpositive factor in front of $v(t,0^+)^2$) $$\tilde{C}((C_{0,-} + C_{1,-} + C_{2,-})\sqrt{\kappa(q)} + (C_{0,-}' + C_{1,-}' + C_{2,-}')) - \dfrac{\omega C_{0,+}}{2}\leq 0,$$ 
\item(to have a nonpositive factor in front of $v(t,0^-)^2$) $$\tilde{C}(\sqrt{\kappa(q)}(C_{1,-} + C_{2,-}) + C_{0,-}' + C_{1,-}' + C_{2,-}') - \dfrac{\omega C_{0,-}}{2} \leq 0,$$ \item(to have a negative factor in front of $E_{-,2}(\zeta_x(t,\cdot))$) $$\tilde{C}|q|C_{1,-} - \dfrac{C_{1,-}'\omega}{2} \leq 0,$$ \item(to have a negative factor in front of $E_{-,1}(v_x(t,\cdot))$)$$\tilde{C}C_{2,-} - \dfrac{C_{1,-}\omega}{2} \leq 0,$$ \item(to have a negative factor in front of $E_{-,2}(\zeta_{xx}(t,\cdot))$)$$\tilde{C}|q|C_{2,-} - \dfrac{C_{2,-}'\omega}{2} \leq 0,$$ \item(to have a nonpositive factor in front of $\zeta(t,0^+)^2$)$$\tilde{C}(C_{0,-}' + C_{1,-}' + C_{2,-}') - \dfrac{\omega C_{0,+}'}{2} \leq 0,$$ \item(to have a nonpositive factor in front of $\zeta_x(t,0^+)^2$)$$\tilde{C}(C_{1,-}' + C_{2,-}') - \dfrac{\omega C_{1,+}'}{2} \leq 0,$$ \item(to have a nonpositive factor in front of $\zeta_{xx}(t,0^-)^2$)$$\tilde{C}C_{2,-}' - \dfrac{\omega C_{2,+}'}{2} \leq 0,$$ \item(to have a nonpositive factor in front of $v_x(t,0^-)^2$)$$\tilde{C}C_{2,-}' - \omega C_{1,-} \leq 0$$ and \item(to have a nonpositive factor in front of $v_x(t,0^+)^2$)$$\tilde{C}C_{2,-} - \omega C_{1,+} \leq 0.$$
\end{itemize} Thus, choosing $C_{0,-}=1=C_{0,+}=C_{1,+}=C_{2,+}=C_{0,+}'=C_{1,+}'=C_{2,+}'$,
  $C_{1,-} = \dfrac{\omega}{4\tilde{C}}$ $C_{2,-} = \min \left(C_{1,-},\dfrac{C_{1,-}\omega}{2\tilde{C}}\right)$
  $C_{0,-}' = \dfrac{\omega}{8\tilde{C}}$ $C_{1,-}' = \min \left(\dfrac{\omega}{8\tilde{C}},\dfrac{\omega C_{1,-}}{\tilde{C}}\right)$ we can obtain the intermediate inequalities as wanted if we choose $q_0$ small enough to make the terms in $|q|$ and $\kappa(q)$ small enough. Thus, for initially smooth and compactly supported functions solutions to our (perturbative) equation, we have that\begin{align*}
  &\sum_{k=0}^2(C_{k,-}\dfrac{dE_{-,1}(\partial^k_xv)}{dt}(t)+C_{k,-}'\dfrac{dE_{-,1}(\partial^k_x\zeta)}{dt}(t) + C_{k,+}\dfrac{dE_+(\partial^k_xv(t))}{dt}(t) + C_{k,+}'\dfrac{dE_+(\partial^k_x\zeta(t)))}{dt}(t) \\ &\leq -\dfrac{\omega}{2}\sum_{k=0}^2(E_{-,1}(\partial^k_xv(t, \cdot))+E_{-,1}(\partial^k_x\zeta(t, \cdot)) + E_+(\partial^k_xv(t, \cdot)) + E_+(\partial^k_x\zeta(t, \cdot))).
  \end{align*}
The desired results are then direct consequences of this inequality.
\end{proof}
Now we can present the proof of the main theorem of this part
\begin{proof}
Take the constants obtained in Proposition \ref{prop:constants}. Given $\eta$, $\eps$, $M$ and $\vartheta$, we will need to apply a continuity argument. First, we work with smooth initial data, compactly supported away from the shock. As the energy is decaying exponentially fast as long as the $W^{1,\infty}$ norm is small enough, and as the energy is a norm on $H^2_{\eps}$ equivalent to the one defined before, we have the result as long as $\|(v,\zeta)\|_{W^{1,\infty}}$ remains small enough. To guarantee that $(v,\zeta)$ has a $W^{1,\infty}$ norm that remains small enough, we use the Sobolev embedding theorem and a continuity argument, shrinking $\delta_0$ if needed. After that, we can just obtain the result for initial data that may not be smooth nor compactly supported (but supported away from the shock) through a density argument, using the continuity of the flow from Theorem \ref{th1}. We note $(v,\zeta)$ the associated maximum solution to a given initial data satisfying the theorem \ref{th1}, as well as $(v_n,\zeta_n)_n$ the maximum solutions associated with a sequence of initial data smooth and with compact supports in $\R^*$ that approximate $(v(0, \cdot),\zeta(0, \cdot))$ in $H^2_{\eps}$. Given that $(v(0, \cdot),\zeta(0, \cdot))$ is small enough in $H^2_{\eps}$, we have that $(v_n,\zeta_n)$ is defined on $\R_+$ for every $n$ big enough, and so we obtain the convergence in $L^2$ of the $(v_n,\zeta_n)$ for every $t$ in the interval of existence of $(v,\zeta)$, and, as we have obtained that $(v_n,\zeta_n)_n$ is bounded in $H^2_{\eps}$ for $n$ big enough, the limit is also bounded in $H^2_{\eps}$, and we obtain that the sequence goes to $0$ with the rate of convergence we were aiming for.
\end{proof}
\begin{proposition}
For any $q$, not necessarily in $[-q_0,q_0]$, such that the wave considered in Proposition \ref{exstwave}, with fixed $f$, $k$ and $\uu(0^-)$, exists then, there exists $\vartheta > 0$, $\tilde{C} > 0$ and $\eta > 0$ such that for every perturbative solutions $(v,\zeta)$ defined on some time interval $[0,T]$ where $T > 0$ such that $\|(v(t, \cdot),\zeta(t, \cdot))\|_{W^{1,\infty}(\R^*)} < \eta$ on $[0,T]$, we have that, for every $t \in [0,T]$ $$\|(v(t, \cdot),\zeta(t, \cdot)\|_{H^2_{\eps}(\R^*)} \leq Ce^{-\vartheta t}\|(v(0),\zeta(0)\|_{H^2_{\eps}(\R^*)} + \int^t_0 C e^{-\vartheta (t-s)}\|(v(s),\zeta(s))\|_{L^2_{\eps}(\R^*)}ds$$
\end{proposition}
The proof follows the same idea as the result before, just this time setting aside the lower order terms, that is the terms of the form $\|v\|_{L^2_{\eps}}^2$ or $\|\zeta\|_{L^2_{\eps}}$.  The main new ingredient in this proof compared to the proof of the stability result is the use of the Sobolev embedding theorem as follows (where $\delta_2$ is a given positive constant that will be adjusted later, and $h$ is an element of $H^1(\R^*)$)
\begin{align*}
|h(0^-)|^2 \leq \delta_2 \|h_x\|_{L^2(\R^-)}^2 + \dfrac{\|h\|_{L^2(\R^-)}^2}{4 \delta_2}, \\
|h(0^+)|^2 \leq \delta_2 \|h_x\|_{L^2(\R^+)}^2 + \dfrac{\|h\|_{L^2(\R^+)}^2}{4 \delta_2}.
\end{align*}
As for the complete proof, it will not be given here. We just focus on how to obtain bounds at the $H^1$ level in the linear case, omitting the proof of the $H^2$ nonlinear problem. \\

We have the following equations for the spatial derivatives of $(v,\zeta)$ (for the linear problem)
\begin{align*}
	(v_x)_t&=(\sigma - f'(\uu))v_{xx} - 2f''(\uu)\uu'v_x - (f^{(3)}(\uu)\uu'^2 + f''(\uu)\uu'')v + kq\zeta_x + (\psi' - \sigma)\uu'', \quad x\in\RR^-,\\
	(\zeta_x)_t&=\sigma \zeta_{xx} + (\psi' - \sigma)\uu'' - k\zeta_x, \quad x\in\RR^-,\\
	(v_x)_t&=(\sigma - f'(0))v_{xx}, \quad x\in\RR^+,\\
	(\zeta_x)_t&=\sigma\zeta_{xx}, \quad x\in\RR^+.
\end{align*}
This gives us the following energy estimates for the functions and their derivatives for smooth solutions (on $[0,T] \times \R^*$) for all $t \in (0,T)$
\begin{align*}
\dfrac{1}{2}\dfrac{dE_+(v)}{dt}(t) &\leq 0,
\end{align*}
\begin{align*}
\dfrac{1}{2}\dfrac{dE_+(\zeta)}{dt}(t) &\leq -\dfrac{\sigma}{2}(\zeta(t,0^+))^2,
\end{align*}
\begin{align*}
\dfrac{1}{2}\dfrac{dE_{-,1}(v)}{dt}(t) &\leq \dfrac{k|q|E_{-,1}(v(t,\cdot))}{2} + \dfrac{k|q|E_{-,1}(\zeta(t,\cdot))}{2} + \dfrac{E_{-,1}(v(t,\cdot))}{2}\\ & \; + \left(\delta_2\|v_x(t,\cdot)\|_{L^2(\R^+)}^2 + \delta_2 \|v_x(t,\cdot)\|_{L^2(\R^-)}^2 + \dfrac{\|v(t,\cdot)\|_{L^2(\R^-)}}{4\delta_2} \right. \\ & \left. + \dfrac{\|v(t,\cdot)\|_{L^2(\R^+)}}{4\delta_2}\right)\dfrac{E_{-,1}(\uu')}{2}+ \dfrac{C_f\sqrt{\kappa(q)}E_{-,1}(v(t,\cdot))}{2},
\end{align*}
\begin{align*}
\dfrac{1}{2}\dfrac{dE_{-,2}(\zeta)}{dt}(t) &\leq \dfrac{C_f}{2}\left(\delta_2\|v_x(t,\cdot)\|_{L^2(\R^+)}^2 + \delta_2\|v_x(t,\cdot)\|_{L^2(\R^-)}^2 + \dfrac{\|v(t,\cdot)\|_{L^2(\R^-)}^2}{4\delta_2} \right. \\ & \left. + \dfrac{\|v(t,\cdot)\|_{L^2(\R^+)}}{4\delta_2}\right) E_{-,2}(\uz') + \dfrac{C_fE_{-,2}(\zeta(t,\cdot))}{2} 
\end{align*}
\begin{align*}
\dfrac{1}{2}\dfrac{dE_+(v_x)}{dt}(t) &\leq -\dfrac{(\sigma - f'(0))}{2}v_x(t,0^+)^2 - \dfrac{(\sigma - f'(0))\eps E_+(v_x(t,\cdot))}{2}
\end{align*}
\begin{align*}
\dfrac{1}{2}\dfrac{dE_+(\zeta_x)}{dt}(t) &\leq -\dfrac{\sigma (\zeta_x(t,0^+))^2}{2} - \dfrac{\sigma \eps E_+(\zeta_x)}{2}
\end{align*}
\begin{align*}
\dfrac{1}{2}\dfrac{dE_{-,2}(\zeta_x)}{dt}(t) &\leq + \dfrac{\sigma (\zeta_x(t,0^-))^2}{2} - \dfrac{2k - \eps \sigma}{2}E_{-,2}(\zeta_x(t,\cdot)) + \dfrac{C_fE_{-,2}( \sqrt{|\uz''|} )}{2} \left(\delta_2\|v_x(t,\cdot)\|_{L^2(\R^+)}^2 \right. \\ & \left. + \delta_2\|v_x(t,\cdot)\|_{L^2(\R^-)}^2 + \dfrac{\|v(t,\cdot)\|_{L^2(\R^-)}^2}{4\delta_2} + \dfrac{\|v(t,\cdot)\|_{L^2(\R^+)}}{4\delta_2} \right)+ \dfrac{C_fE_{-,2}(\sqrt{|\uz''|}\zeta_x)}{2} \\ & + \int_{\R^-}\dfrac{-C\sigma e^{\eps \cdot}(\zeta_x)^2}{2} \varrho_{-,2},
\end{align*}
\begin{align*}
\dfrac{1}{2}\dfrac{dE_{-,1}(v_x)}{dt}(t) &\leq \dfrac{(\sigma - f'(u_0))}{2}(v_x(t,0^-))^2 - \dfrac{-\eps (\sigma - f'(u_0))}{2}E_{-,1}(v_x(t,0^-))^2 \\ & - \int_{\R^-}\dfrac{Ce^{\eps \cdot}(v_x(t,\cdot))^2(f'(u_0)-\sigma)}{2}\varrho_{-,1} + 2 C_f\int_{\R^-}|\uu'|v_x(t,\cdot)^2\varrho_{1,-} \\ & + \dfrac{C_f}{2} \int_{\R^-}(\uu'^2 + |\uu''|)v_x(t,\cdot)^2\varrho_{-,1} + \dfrac{C_f}{2} \int_{\R^-}v(t,\cdot)^2\varrho_{-,1} \\ & + \dfrac{k|q|E_{-,1}(\zeta_x(t,\cdot))}{2\delta_1} + \dfrac{k|q|\delta_1E_{-,1}(v_x(t,\cdot))}{2} + \left(\delta_2\|v_x(t,\cdot)\|_{L^2(\R^+)}^2 \right. \\ & \left. + \delta_2\|v_x(t,\cdot)\|_{L^2(\R^-)}^2 + \dfrac{\|v(t,\cdot)\|_{L^2(\R^-)}^2}{4\delta_2} + \dfrac{\|v(t,\cdot)\|_{L^2(\R^+)}}{4\delta_2}\right) \dfrac{E_{-,1}(\sqrt{|\uu''|})}{2} \\ & + \dfrac{E_{-,1}(\sqrt{|\uu''|}v_x(t,\cdot))}{2}.
\end{align*}
Thus, choosing $C$, $\delta_1$ and $\delta_2$ as $C=\max\left(\dfrac{C_fk^2}{\sigma^3},\dfrac{\kappa(q) + 2C_f(3\kappa(q) + \kappa(q)^2)}{f'(u_0) - \sigma}\right)$ and $0 < \delta_2 < \min \left( \dfrac{\eps(f'(u_0) - \sigma)}{8E_{-,1}(\sqrt{|\uu''|})},\dfrac{\eps (f'(u_0) - \sigma)(2k - \eps\sigma)\eps \sigma}{16C_fE_{-,2}(\sqrt{|\uz''|})k^2q^2e^{2C\eps^{-1}}} \right) $ as well as $\delta_1=\dfrac{\eps (f'(u_0) - \sigma)}{4k|q|}$ (it does not work for $q=0$, but this case is contained in the stability result proved earlier, so we assume that $q \neq 0$) we will be able to conclude the desired damping estimate. For any $t \in (0,T)$ we obtain \begin{align*}
\dfrac{k|q|e^{2C\eps^{-1}}}{(2k-\eps\sigma)\delta_1}\dfrac{d(E_{-,2}(\zeta_x))}{dt}(t) + \dfrac{1}{2}\dfrac{d(E_{-,1}(v_x))}{dt}(t) \leq& \dfrac{-\eps (f'(u_0) - \sigma)E_{-,1}(v_x(t,\cdot))}{4} - \dfrac{k|q|e^{2C\eps^{-1}}}{4\delta_1}E_{-,2}(\zeta_x(t,\cdot)) \\ & + \dfrac{C_fE_{-,1}(v(t,\cdot))}{2} + \dfrac{\|v(t,\cdot)\|_{L^2(\R)}^2E_{-,1}(\sqrt{|\uu''|})}{8\delta_2} \\ & + \dfrac{k|q|e^{2C\eps^{-1}}\|v(t,\cdot)\|_{L^2(\R)}^2E_{-,1}(\sqrt{|\uz''|})}{4\delta_2(2k-\eps\sigma)\delta_1} \\ & + \left(\|v_x(t,\cdot)\|_{L^2(\R^+)}^2\delta_2 + \dfrac{\|v(t,\cdot)\|_{L^2(\R)}^2}{4\delta_2}\right) \times \\ &\left( \dfrac{C_f(E_{-,1}(\sqrt{|\uu''|}) + E_{-,2}(\sqrt{|\uz''|})}{2} \right) \\ &+ \dfrac{k|q|e^{2C\eps^{-1}}(\zeta_x(t,0^+))^2\sigma}{2\delta_1(2k - \eps \sigma)}.
\end{align*} Hence, by now taking the full energy estimates (we write $\mathcal{E}(v,\zeta):=E_{-,1}(v) + E_{-,2}(\zeta) + E_{-,1}(v) + \dfrac{E_{-,2}(\zeta)k|q|e^{2C\eps^{-1}}}{k-\dfrac{\eps\sigma}{2}\delta_1} + E_+(v) + ME_+(\zeta) + E_+(v_x) + ME_+(\zeta_x)$, where $M$ is a positive constant to be determined, the following bounds are obtained
\begin{align*}
\dfrac{1}{2}\dfrac{d\mathcal{E}(v,\zeta)}{dt}(t) \leq& -\dfrac{M\sigma}{2}((\zeta(t,0^+))^2 + (\zeta_x(t,0^+))^2) - \dfrac{\eps (f'(u_0) - \sigma)E_{-,1}(v_x(t,\cdot))}{4} - \dfrac{k|q|e^{2C\eps^{-1}}}{4\delta_1}E_{-,2}(\zeta_x(t,\cdot)) \\ & + \left(\|v_x(t,\cdot)\|_{L^2(\R^+)}^2\delta_2\right) \times \left( \dfrac{C_f(E_{-,1}(\sqrt{|\uu''|}) + E_{-,2}(\sqrt{|\uz''|})}{2} \right) + \dfrac{k|q|e^{2C\eps^{-1}}(\zeta_x(t,0^+))^2\sigma}{2\delta_1(2k - \eps \sigma)} \\ & + \dfrac{\sigma (\zeta(t,0^+))^2}{2} - \dfrac{\eps (\sigma - f'(0))}{2}E_+(v_x(t,\cdot)) + \tilde{C}(\|v(t,\cdot)\|_{L^2(\R)}^2 + \|\zeta(t,\cdot)\|_{L^2(\R)}^2) \\ & - \dfrac{M\sigma \eps E_+(\zeta_x(t,\cdot))}{2}.
\end{align*}
For some constant $\tilde{C}$ which depends on $\delta_2$. Thus, for $M$ big enough and $\delta_2$ small enough, we obtain \begin{align*}
\dfrac{1}{2}\dfrac{d(\mathcal{E}(v,\zeta))}{dt}(t) \leq \dfrac{\gamma\mathcal{E}(v,\zeta)}{dt}(t) + \tilde{C}(\|v(t,\cdot)\|_{L^2(\R)}^2 + \|\zeta(t,\cdot)\|_{L^2(\R)}^2)
\end{align*}
for some $\gamma<0$.
\subsection{Negative propagation speed}
We now examine the case $\sigma<0$.
\begin{proposition}
Let $\alpha \in \R$, such that $\alpha \neq 0$ or for all neighborhoods $V$ of $0$, $f'$ is not constant on $V$. We assume that $\sigma < 0$. There exists $\eps_0 > 0$ small enough such that there exists $(v_{0,n},\zeta_{0,n})_{n \in \N} \in C^{\infty}_c(\R)^{\N}$ that goes to $0$ in $H^2_{\alpha}$, supported away from $0$, and such that the induced solution remains supported away from the shock on some time interval $[0,T_{n,\alpha})$ and  $\sup_{t \in [0,T_{n,\alpha})}\|(v_n(t,\cdot),\zeta_n(t,\cdot)\|_{H^2_{\alpha}(\R)}$ is larger than $\eps_0$.
\end{proposition}
\begin{proof}
Case $\alpha \leq 0$: \\
When $v_t+(f(v))_x=0$ has a compactly supported initial perturbation, we can proceed as follows. We consider $(\phi_n)_n \in H^2(\R)$ fixed and all nonzero, all initial data such that the induced solutions blow up in finite time in $W^{1,\infty}$, and with the sequence of their $H^2$ norms that goes to $0$ when $n$ goes to $+\infty$. We also assume that their support is in $[0,1]$. \\

We consider the sequence of initial perturbation $(v_n,0)$ where $v_n : x \mapsto \phi_n(\cdot-a_n)$, where $(a_n)_n$ is such that the solution induced by $\phi_n$ has its support included in $[0,T_n) \times (-a_n,+\infty)$ (where $T_n$ is its blow up time). We are now reduced to the case of a scalar conservation law. For $n$ big enough, we have that the solution blows up in finite time in $W^{1,\infty}$, and, thus, in $H^2_{\alpha}$ (as the support of the perturbation at time remains in a compact set $K_n$ independent of $t$, we can apply the the Sobolev embedding to obtain a bound on the unweighted $W^{1,\infty}$ norm).  \\

Now, we assume that $f'$ is constant, and $\alpha<0$. By assumption, we have that $f'(0)<\sigma$. Thus, by taking any initial data nonzero localized on the right of the shock, we have that the solutions grows in $L^2_{\alpha}$ exponentially fast in time as long as the solution is localized on the right of the shock. In particular, by translating some bump function far enough to the right, we have the desired result. \\

The initial data given by $(\ul{u}+v_n,\ul{z}) $(when taking $n$ big enough, after multiplying it by some small $\delta>0$) either induces a solution that blows up in finite time in $W^{1,\infty}$, or a solution that remains small in $L^{\infty}$, is defined globally.\\
 
Case $\alpha > 0$:  \\ This time, we will not obtain explosion in finite time. However, with the initial data $( 0 , (n+1)^{-1}\phi )_n$, for $t$ such that the solution continues to exist on $[0,t]$ we have that: $\zeta_n(t,x)=\dfrac{\phi(x-t\sigma)}{n+1}$ for all $x>0$. Furthermore, as $\zeta_n(t,x)=0$ for all $x<0$ and $v_n(t,x)=0$ for all $x\in \R^*$, we obtain the desired instability result: the solution is defined globally, but its $\|\cdot\|_{L^2_{\alpha}}$ norm goes to $+\infty$ as $t$ goes to $+\infty$.
\end{proof}
\section{Singularity formation for ZND}\label{s:znd}
In this section we prove the negative result Theorem \ref{thm:generalblowup} and Corollary \ref{cor:ZNDBlowup} as well as a similar blowup result for \emph{un}weighted norms for the Majda model. We will also discuss the situation with weighted norms for ZND.
\subsection{Setup}
Let $A(u):\RR^n\to M_n(\RR)$ be a smooth matrix function with the property that there exists a $\delta>0$ so that $A(u)$ is strictly hyperbolic for $|u|\leq \delta$. Further assume that the eigenvalues $\L_i(u)$ of $A(u)$ may be ordered as
\be\label{eq:AstrictHyperbolic}
\L_n(u)<\L_{n-1}(u)<...<\L_1(u),
\ee
and that each $\L_j$ is simple. By shrinking $\delta$ if necessary, we may further assume that $\L_i(v)<\L_j(u)$ if $j<i$ and $|u|,|v|\leq\delta$. Assume further that there exists a smooth function $F:\RR^n\to\RR^n$ with
\begin{equation*}
	A(U)=D_UF(U),
\end{equation*}
so that
\be
	U_t+F(U)_x=U_t+A(U)U_x=0,
\ee
is a system of conservation laws.\\

Let $\underline{U}(x-\s t)$ be a shock solution to
\be\label{eq:unperturbedeqn}
U_t+A(U)U_x=0,
\ee
for some constant $\s$ and such that $||\underline{U}||_{\infty}<\delta$ and $\underline{U}$ and all of its derivatives are exponentially decaying in space. For our purposes, we need $A$ to be strictly hyperbolic with at least one genuinely nonlinear field. As both the required conditions on $A$ are invariant under Galilean changes of coordinates, we may without essential loss of generality assume that the underlying shock speed $\s=0$. Suppose we have a solution to \eqref{eq:unperturbedeqn} of the form 
\be
U=\underline{U}+\hat{U},
\ee
with $\hat{U}$ small. Then the evolution for the perturbed solution $U$ is given by
\be\label{eq:perturbedderiv1}
u_t+A(U)U_x=0.
\ee
Subtracting the equation for the shock
\be\label{eq:shockeqn}
A(\underline{U})\underline{U}_x=0,
\ee
from \eqref{eq:perturbedderiv1} produces
\be\label{eq:perturbedderiv2}
\hat{U}_t+(A(\underline{U}+\hat{U})-A(\underline{U}))\underline{U}_x+A(\underline{U}+\hat{U})\hat{U}_x=0.
\ee
By applying the fundamental theorem of calculus twice, we see that
\ba\label{eq:perturbedderiv3FTOC}
A(\underline{U}+\hat{U})-A(\underline{U})&=\int_0^1 D_uA(\underline{U}+s\hat{U})\hat{U}ds=\tilde{G}(x,\hat{U})\hat{U},\\
A(\underline{U}+\hat{U})&=A(\hat{U})+\int_0^1D_uA(\hat{U}+s\underline{U})\underline{U}ds=A(\hat{U})+B(x,\hat{U}),
\ea
where $D_u$ denotes the Fr\'echet derivative with respect to $u$, $\tilde{G}$ is a bilinear form and $B$ is a matrix. We have that $\tilde{G}$ and $B$ are piecewise smooth, smooth for $x$ large, and that $B$ and all of its derivatives are exponentially decaying in $x$.
\be\label{eq:perturbedeqn}
\hat{U}_t+(A(\hat{U})+B(x,\hat{U}))\hat{U}_x+G(x,\hat{U})\hat{U}=0,
\ee
for $G(x,\hat{U})\hat{U}=\tilde{G}(x,\hat{U})(\hat{U},\underline{U}_x)$. Note that $G$ and all of its derivatives are exponentially decaying in $x$ due to the presence of $\underline{U}_x$. Now that we have the equation for the perturbation \eqref{eq:perturbedeqn}, we are now going to call the perturbation $u$ in order to match the notation of \cite{J} more closely.\\

We note that there exists an $\e>0$ and $R>0$ depending on $B(x,u)$ and $G(x,u)$ so that $B$ and $G$ are smooth for $|x|\geq R$ and such that for all $|x|\geq R$ one has that
\be\label{eq:BGsmall}
|B(x,u)|, |G(x,u)|<\e,
\ee
and the matrix
\be
\cA(x,u):=A(u)+B(x,u),
\ee
is strictly hyperbolic with all simple eigenvalues $\l_j(x,u)$ which smooth depend on $x,u$ \cite{K} and satisfy $\l_i(x,v)<\l_j(x,u)$ for $|u|,|v|\leq\delta$, shrinking $\delta$ slightly if necessary while preserving $||\underline{u}||_{\infty}<\delta$, and $j<i$. This is due to the exponential decay of $B,G$ with respect to $x$.
From \cite{K}, we know that we can also find left and right eigenvectors $\eta_j(x,u)$ and $\xi^j(x,u)$ of $\cA(x,u)$ associated to the eigenvalues $\l_j(x,u)$ with the biorthogonality condition
\be\label{eq:biorthogonality}
\eta_i(x,u)\xi^j(x,u)=\delta_{ij}.
\ee
Further assume the normalization condition on the $\eta_i$
\be\label{eq:etanormalization}
\eta_i\cdot\eta_i=1.
\ee
If a solution $\hat{U}(x,t)$ to \eqref{eq:perturbedeqn} remains supported in $\{x: |x|\geq 2R \}$ for all time $t$ on some sufficiently long time interval $0\leq t\leq T$, then we may smooth out $B$ and $G$ in a such a way as to ensure that \eqref{eq:BGsmall} holds everywhere in $x$ by choosing smooth approximations $\hat{B}$ and $\hat{G}$ such that $\supp(B-\hat{B} ),\supp(G-\hat{G})\subset\{x:|x|\leq\frac{3}{2}R \}\times\{u: |u|\leq\delta \}$. As we're only trying to show blow up of a specific solution which will be supported far from $x=0$ on the desired time of existence, there is no real harm in assuming that $B$ and $G$ are smooth everywhere and that the inequality \eqref{eq:BGsmall} holds for all $x$. This is done purely to avoid some technicalities arising when trying to solve for the characteristics of \eqref{eq:perturbedeqn} near $x=0$. \\

We introduce coefficients $b_{ij}(x,u)$ and $c_{ijk}(x,u)$ in order to describe the gradient of $\cA(x,u)$. In particular, we define
\be\label{def:bij}
b_{ij}(x,u):=\eta_i(x,u)\frac{\d \cA}{\d x}(x,u)\xi^j(x,u),
\ee
and
\be\label{def:cijk}
c_{ijk}(x,u):=\eta_i(x,u)\Big(\frac{d}{d s}\cA(x,u+s\xi^k)|_{s=0}\Big)\xi^j(x,u).
\ee
\begin{remark}
	Morally, $c_{ijk}(x,u)$ can be thought of as
	\be
	c_{ijk}(x,u)=\eta_i(x,u)\frac{\d\cA}{\d u_k}(x,u)\xi^j(x,u),
	\ee
	when one writes $u$ as
	\be
	u=\sum_{k=1}^nu_k\xi^k.
	\ee
\end{remark}
From this, we get the following identities
\ba\label{eq:KatoFormulasLambdaEta}
d\l_i(x,u)&=b_{ii}(x,u)dx+\sum_{k=1}^n c_{iik}(x,u)(du)_k,\\
d\eta_i(x,u)&=\sum_{k\not=i}\frac{b_{ik}(x,u)}{\l_i(x,u)-\l_k(x,u)}( (\eta_k\cdot\eta_i)\eta_i-\eta_k ) dx+\\
&\quad+\sum_{k,m \atop k\not=i}\frac{1}{\l_k(x,u)-\l_i(x,u)}c_{ikm}(x,u)(du)_m( (\eta_k\cdot\eta_i)\eta_i-\eta_k ),
\ea
where we've written $du$ as
\be
du=\sum_{k=1}^n(du)_k\xi^k=\sum_{k=1}^n(\eta_k(x,u)du )\xi^k(x,u).
\ee
We record the decay in $B$ and $G$ as
\be\label{eq:expdecayrate}
|\d_x^\a \d_u^\b B(x,u)|, |\d_x^\a \d_u^\b G(x,u)|\leq C_{\a\b} e^{-c|x|},
\ee
for $|x|$ bounded away from 0.
\subsection{Characteristics}
For each $i=1,...,n$ we let $X_i(x,t)$ be the solution to the ODE
\ba\label{eq:ithCharacteristicEqn}
\frac{\d X_i}{\d t}(x,t)&=\l_i(X_i(x,t),u(X_i(x,t),t)),\\
X_i(x,0)&=x.
\ea
Let $\overline{\l}$ be defined as
\be\label{eq:maxspeed}
\overline{\l}:=\max_{1\leq i \leq n}\sup_{x\in \RR \atop |u|\leq \delta}|\l_i(x,u)|<\infty.
\ee
The observation that $\overline{\l}<\infty$ follows from $||\cA||_{\infty}<\infty$. Note that $\overline{\l}<\infty$ holds even when $B$ arises from linearizing about a shock as $||B||_{\infty}<\infty$. Define parameters $\nu_i,\mu_i$ by
\ba
\nu_i&:=\inf_{x\in\RR,\ |u|\leq \delta} \l_i(x,u),\\
\mu_i&:=\sup_{x\in\RR,\ |u|\leq \delta} \l_i(x,u).
\ea
Further define $\s$ by
\be\label{def:sigma}
\s:=\min_{k<i}\nu_k-\mu_i>0.
\ee
Let $u(x,t)$ be a solution to
\be
u_t+\cA(x,u)u_x+G(x,u)u=0
\ee
such that for some $\infty>T>0$ one has the bound
\be
\sup_{0\leq t\leq T\atop x\in\RR} |u(x,t)|\leq \delta.
\ee
Thus if $u(x,0)$ is compactly supported in $[\a_0,\b_0]$, we then have by the method of characteristics that the solution $u(x,t)$ is compactly supported in $[\a_0-\overline{\l}T,\b_0+\overline{\l}T]$. As such, if $[\a_0,\b_0]$ is sufficiently far from $x=0$ then $u(x,t)$ remains supported in $|x|\geq 2R$ for all $t$ and hence we may smooth out $B$ and $G$. \\

Let $w:=u_x$ and write $w$ as
\be
w=\sum_{i=1}^n w_i\xi^i.
\ee
From this, we can write the evolution of $u$ along the $i$-th characteristic as
\be\label{eq:ualongichar}
\frac{d u}{dt}=u_t+\l_i u_x=(\l_i-\cA(x,u))u_x-G(x,u)u=\sum_{k=1}^n(\l_i-\l_k)w_k\xi^k-G(x,u)u.
\ee
Correspondingly, we have the evolution of $\l_i$ along the $i$-th characteristic given by
\ba\label{eq:lambdaialongichar}
\frac{d\l_i}{dt}&=b_{ii}\frac{\d X_i}{\d t}+\sum_{k=1}^n c_{iik}\eta_k\frac{du}{dt},\\
&\quad=\l_ib_{ii}+\sum_{k=1}^n \big( (\l_i-\l_k)w_k-\eta_k G(x,u)u\big)c_{iik}.
\ea
Similarly, the evolution of $\eta_i$ along the $i$-th characteristic is given by
\ba\label{eq:etaialongichar}
\frac{d \eta_i}{dt}&=\sum_{k\not=i}\frac{b_{ik}}{\l_i-\l_k}( (\eta_k\cdot\eta_i)\eta_i-\eta_k )\frac{\d X_i}{\d t}+\sum_{k,m \atop k\not=i}\frac{1}{\l_k-\l_i}c_{ikm}(\eta_m\frac{du}{dt})( (\eta_k\cdot\eta_i)\eta_i-\eta_k ),\\
&\quad=\sum_{k\not=i}\frac{\l_i }{\l_i-\l_k}b_{ik}( (\eta_k\cdot\eta_i)\eta_i-\eta_k )+\sum_{k,m\atop k\not=i}\Bigg(\frac{(\l_i-\l_m)}{\l_k-\l_i}c_{ikm}w_m ( (\eta_k\cdot\eta_i)\eta_i-\eta_k )-\\
&\quad-\frac{c_{ikm}}{\l_k-\l_i}\eta_m G(x,u)u( (\eta_k\cdot\eta_i)\eta_i-\eta_k )\Bigg).
\ea
This leads us to the evolution of $w_i$ along the $i$-th characteristic as
\ba
\frac{d w_i}{dt}&=\eta_i(u_{xt}+\l_iu_{xx})+\frac{d\eta_i}{dt}u_x\\
&\quad=\eta_i((\l_i-\cA)w_x-\cA_xw)+\frac{d \eta_i}{dt}w-\eta_i(G(x,u)u )_x.
\ea
As $\eta_i$ was assumed to be a left eigenvector of $\cA$ with eigenvalue $\l_i$, we are left with
\be
\frac{d w_i}{dt}=\sum_{k=1}^n(\frac{d\eta_i}{dt}\xi^k )w_k-\sum_{k=1}^n\eta_i\cA_x\xi^kw_k-\eta_i( G_x(x,u)u+G_U(x,u)(u,w)+G(x,u)w).
\ee
Plugging in \eqref{eq:etaialongichar}, we find that
\ba
\frac{d w_i}{dt}&=\sum_{k\not=i}\frac{b_{ik}\l_i}{\l_k-\l_i}w_k+\sum_{k,m\atop k\not= i}\frac{\l_i-\l_m}{\l_k-\l_i}c_{ikm}((\eta_i\cdot\eta_k)w_mw_i-w_kw_m)=\\
&\quad-\sum_{k=1}^n (\eta_i\cA_x\xi^k+\eta_i G\xi^k+\eta_i G_U(u,\xi^k)) w_k-\\
&\quad- \sum_{k,m \atop k\not =i}\frac{c_{ikm}}{\l_k-\l_i}\eta_m G(x,u)u ((\eta_i\cdot\eta_k) w_i-w_k)-\eta_i G_x(x,u)u .
\ea
Now plugging in the expressions from \eqref{def:bij} and \eqref{def:cijk}, we get
\ba
\frac{d w_i}{dt}&=\sum_{k\not=i}\frac{b_{ik}\l_i}{\l_k-\l_i}w_k+\sum_{k,m\atop k\not= i}\frac{\l_i-\l_m}{\l_k-\l_i}c_{ikm}((\eta_i\cdot\eta_k)w_mw_i-w_kw_m)-\\
&\quad-\sum_{k=1}^n b_{ik}w_k-\sum_{k,m}c_{ikm}w_kw_m-\sum_{k}(\eta_i G\xi^k+\eta_i G_U(u,\xi^k) )w_k-\\
&\quad-\sum_{k,m \atop k\not =i}\frac{c_{ikm}}{\l_k-\l_i}\eta_m G(x,u)u ((\eta_i\cdot\eta_k) w_i-w_k)-\eta_i G_x(x,u)u.
\ea
Introducing coefficients $\g_{ikm}$ as in \cite{J}, and new coefficients $\zeta_{ik}$ and $\kappa_i$ we get
\be\label{eq:wichareqn}
\frac{d w_i}{dt}=\sum_{k=1}^n \zeta_{ik} w_k+\sum_{k,m}\g_{ikm}w_kw_m+\kappa_i u.
\ee
We note that each $\zeta_{ik}$ and $\kappa_i$ satisfies
\be
|\d_x^\a\d_u^\b\zeta_{ik}(x,u)|,|\d_x^\a \d_u^\b \kappa_i(x,u)|\leq C_{\a\b}e^{-c|x|}
\ee
for $c>0$ the same constant as in \eqref{eq:expdecayrate}, coming from the exponential decay of $B$ and $G$.
\begin{rem}
	Comparing with \cite{B}, our characteristic equation \eqref{eq:wichareqn} takes a very similar form with the key difference being that here we have a forcing term independent of $w$. The forcing term arises due to the spatial inhomogeneity induced by linearizing \eqref{eq:unperturbedeqn} about a fixed non-constant solution.
\end{rem}
For initial data which is compactly supported in an interval $[\a_0,\b_0]$, we let $\cR_i$ denote the region foliated by the characteristics, i.e.
\be
\cR_i=\{(X_i(x,t),t): x\in[\a_0,\b_0], 0\leq t\leq T \}.
\ee
We also define extremal characteristics
\ba
\a_i(t)&=X_i(\a_0,t),\\
\b_i(t)&=X_i(\b_0,t).
\ea
As in the unperturbed case, when $\s$ defined in \eqref{def:sigma} is positive, there exists a $t_0>0$ so that $\cR_i\cap\cR_j\subset\{(x,t): 0\leq t\leq t_0 \}$. Moreover, this $t_0$ is comparable to the initial width of the support $s_0:=\b_0-\a_0$.
\subsection{Estimates}
We let $u(x,t)$ be a $C^2$ solution to \eqref{eq:perturbedeqn}, which a priori remains bounded by $\delta$ for all $0\leq t\leq T$ and all $x$. Assume that the initial data is small in the sense that
\be
\theta:=\sup_x s_0^2|u_{xx}(x,0)|\ll 1.
\ee
In particular, assume that $\theta\leq \frac{1}{2}\delta$. We may further assume that there exists an index $i$ and a $y\in[\a_0,\b_0]$ so that $w_i(y,0)=\max_i\sup_x |w_i(x,0)|=W_0$. A simple consequence of the fundamental theorem of calculus is that
\be\label{eq:John48}
s_0W_0=\cO(\theta).
\ee
In this section, $X=\cO(Y)$ will mean that $X\leq CY$ for some constant $C$ that only depends on the matrices $A(u)$, $B(x,u)$, $G(x,u)$, $\delta$ and $\e$ for $\theta$ sufficiently small.
\begin{remark}
	The argument in John does not make the assumption that $W_0$ is achieved by some $i$ and $y$ with $w_i(y,0)>0$ in order to show that the blowup is generic. As the initial data we will show blowup for will not be generic due to being spatially supported far from $x=0$, we are free to assume that $W_0$ is achieved when $w_i$ is positive. As another remark, the coefficient $\g_{iii}$ that leads to Riccati-type blowup can be arranged to be positive by changing the sign of $\eta_i$, $\xi^i$. If one assumes that $\g_{iii}$ is negative, then one assumes that there is an $i$ and a $y$ such that $-w_i(y,0)=\max_i \sup_x|w_i(x,0)|$.
\end{remark}
The overall plan will be to follow the argument in \cite{J} as closely as possible while keeping track of the dependence on the constants on the distance of the support of $u(x,0)$ to 0 and the time of existence $T$. We are going to allow ourselves the freedom to choose the distance of the support from 0 on parameters associated to the initial data, such as $W_0$, $s_0$, and ultimately $T$ as well. What we are not allowed to do is allow the distance to depend on information about the solution for $t>0$. We will then construct a specific initial data, originally supported on $[-\frac{1}{2},\frac{1}{2}]$, and then scale it by $\theta$ and translate it out far enough so that all of the necessary distance conditions are satisfied.\\

Before we begin, we let $\cE(d,T)$ be some constant which only depends on the distance $d$ of the support of $u(x,0)$ from 0 and the desired time of existence $T$ such that $\cE(d,T)$ controls all quantities derived from $B(x,u)$ and a sufficient number of their derivatives uniformly on $\supp(u(x,0))+[-\overline{\l}T,\overline{\l}T]$ and in $L^1$. In a similar fashion, we let $\cG(d,T)$ be a constant which controls all quantities derived from $G(x,u)$ and a sufficient number of derivatives uniformly on $\supp(u(x,0))+[-\overline{\l}T,\overline{\l}T]$ and in $L^1$. That is, there is a positive integer $k$ and a large constant $C>0$ depending only on $\cA(x,u)$, $\cG(x,u)$, $\delta$ and $\e$ so that
\ba
	||B||_{W^{k,\infty}}+\sup_{|u|\leq \delta} ||B(\cdot,u)||_{W^{k,1}( (-\infty,\b_0+\overline{\l}T) )}&\leq\frac{1}{C}\cE(d,T),\\
	||G||_{W^{k,\infty}}+\sup_{|u|\leq\delta} ||G(\cdot,u)||_{W^{k,1}(-\infty,\b_0+\overline{\l}T)}&\leq\frac{1}{C}\cG(d,T),
\ea
if $[\a_0,\b_0]\subset(-\infty,0)$ and a similar convention if $[\a_0,\b_0]\subset (0,\infty)$.
We remark that due to the exponential decay in $\underline{u}$, we have the bounds
\be
\cE(d,T),\cG(d,T)\leq Ce^{-cd+cT},
\ee 
for $C,c$ both finite positive constants.\\

Define constants
\begin{subequations}
	\begin{align}
	W&=\sup_i\sup_{(x,t)\atop 0\leq t \leq T} |w_i(x,t)|,\\
	V&=\sup_i \sup_{(x,t)\not\in\cR_i\atop 0\leq t\leq T} |w_i(x,t)|,\\
	U&=\sup_{(x,t)\atop 0\leq t \leq T}|u(x,t)|,\\
	S&=\sup_i \sup_{0\leq t\leq T} (\b_i(t)-\a_i(t)),\\
	J&=\sup_i\sup_{0\leq t\leq T}\int_{\a_i(t)}^{\b_i(t)}|w_i(x,t)|dx.
	\end{align}
\end{subequations}
\begin{lemma}\label{lem:USbounds}
	If $d:=\text{dist}([\a_0,\b_0],0)$ is the distance of the support of $u(x,0)$ from 0, then one can bound $U$ and $S$ by
	\ba\label{eq:keyUSest}
	U&=\cO(s_0V+TV+J), \\
	S&=\cO(s_0+\cE(d,T)+VTS+JT),
	\ea
	Provided that $d\gtrsim\max\{1,T,|\log(s_0)|,|\log(W_0)| \}$, we may assume that \be
	\cE(d,T)\leq\min\{1,s_0, \tilde{c}W_0 \}
	\ee for some $0<\tilde{c}\ll1$. 
\end{lemma}
\begin{proof}
	We recall the proof of the estimate on $U$ from \cite{J}. To bound $u$, one starts with the observation that $u(x,t)$ is supported on $[\a_n(t),\b_1(t)]$ for each $t$. Then applying the fundamental theorem of calculus gives
	\be
	u(x,t)=\int_{\a_n(t)}^x\sum_{k=1}^nw_k(y,t)\xi^k(y,t)dy,
	\ee
	and then bounding $u(x,t)$ by
	\be
	|u(x,t)|\leq \max_k\sup_{(x,u)\atop |u|\leq \delta}|\xi^k(x,u)|\sum_{i=1}^n\int_{\a_n(t)}^{\b_1(t)}|w_i(x,t)|dx.
	\ee
	For each $i$ we may write this as
	\be
	|u(x,t)|\leq C\sum_{i=1}^n\Big(\int_{\a_n(t)}^{\a_i(t)}|w_i(x,t)|dx+\int_{\a_i(t)}^{\b_i(t)}|w_i(x,t)|dx+\int_{\b_i(t)}^{\b_1(t)}|w_i(x,t)|dx \Big).
	\ee
	The middle integral may be bounded by $J$ by definition, and the other two integrals may be bounded by $(\a_i(t)-\a_n(t))V+(\b_1(t)-\b_i(t))V$ in a similar manner. But \be
	|(\a_i(t)-\a_n(t))+(\b_1(t)-\b_i(t))|=|\b_1(t)-\a_n(t)-(\b_i(t)-\a_i(t))|\leq 2|\b_1(t)-\a_n(t)|,
	\ee
	leading to
	\be
	|u(x,t)|\leq 2|\b_1(t)-\a_n(t)|V+J.
	\ee
	To finish the estimate on $U$, we note that
	\be
	\Big|\frac{d (\b_1(t)-\a_n(t))}{dt}\Big|=\big|\l_1(\b_1(t),u(\b_1(t),t))-\l_n(\a_n(t),u(\a_n(t),t))\big|\leq 2\overline{\l},
	\ee
	which gives the estimate $|\b_1(t)-\a_n(t)|\leq s_0+2\overline{\l}T$.\\
	
	The estimate for $S$ works slightly differently from the original argument in \cite{J}. We start in the same way by looking at
	\be\label{eq:Sest1}
	\frac{d (\b_i(t)-\a_i(t))}{dt}=\l_i(\b_i(t),u(\b_i(t),t))-\l_i(\a_i(t),u(\a_i(t),t)).
	\ee
	Applying the fundamental theorem of calculus gives
	\be\label{eq:Sest2}
	\frac{d (\b_i(t)-\a_i(t))}{dt}=\int_{\a_i(t)}^{\b_i(t)}\frac{\d }{\d x}\l_i(x,u(x,t))dx.
	\ee
	Recalling the expansion for $d\l_i$ from \eqref{eq:KatoFormulasLambdaEta}, we find that
	\be\label{eq:Sest3}
	\frac{d(\b_i(t)-\a_i(t))}{dt}=\int_{\a_i(t)}^{\b_i(t)}b_{ii}(x,u(x,t))+\sum_{k=1}^n c_{iim}(x,u(x,t))w_m(x,t) dx.
	\ee
	First, we bound
	\be\label{eq:Sest4}
	\Big|\int_{\a_i(t)}^{\b_i(t)}b_{ii}(x,u(x,t))dx\Big|\leq C\int_{\a_i(t)}^{\b_i(t)}e^{-c|x|}dx\leq C\big(e^{-c|\b_i(t)|}-e^{-c|\a_i(t)|}\big)\leq \cE(d,T),
	\ee
	Choosing $\a_0,\b_0$ sufficiently far from 0, depending on $T$ and $s_0$, one may ensure that
	\be\label{eq:Sest5}
	\int_0^t\Big|\int_{\a_i(s)}^{\b_i(s)}b_{ii}(x,u(x,s))dx\Big|ds\leq\cE(d,T) \leq \min\{s_0,1\}
	\ee
	for all $0\leq t\leq T$.\\
	Turning to the sum in \eqref{eq:Sest3}, we split into two cases depending on whether $t\leq t_0$ or $t\geq t_0$ holds. We start with the latter case, where we see that
	\be\label{eq:Sest6}
	\Big|\int_{\a_i(t)}^{\b_i(t)} \sum_{m=1}^nc_{iim}(x,u)w_m(x,t)dx\Big|\leq \max_{i,m}\sup_{(x,u)\atop |u|\leq \delta}|c_{iim}(x,u)|\sum_{m=1}^n\int_{\a_i(t)}^{\b_i(t)}|w_m(x,t)|dx.
	\ee
	To use the assumption $t\geq t_0$, we recall that that $t\geq t_0$ implies that $(x,t)$ can be in at most one $\cR_k$. As the integral in \eqref{eq:Sest6} is over a horizontal slice of one such $\cR_k$, we can further bound \eqref{eq:Sest6} by
	\be\label{eq:Sest7}
	(n-1)\int_{\a_i(t)}^{\b_i(t)}Vdx+\int_{\a_i(t)}^{\b_i(t)}|w_i(x,t)|dx\leq C(SV+J).
	\ee
	To handle $t\leq t_0$, we first show $w_i$ remains under control for $t\leq t_0$. To do this, we first introduce auxiliary parameters $Z$ and $\Gamma$ as follows
	\ba\label{def:gronwallparam}
	K(d,T)&=\sup_{x\in [\a_0-\overline{\l}T,\b_0+\overline{\l}T]\atop |u|\leq\delta } \sum_{i=1}^n|\kappa_i(x,u)|,\\
	Z&=\sup_{(x,u)\atop|u|\leq \delta}\sum_{ik}|\zeta_{ik}(x,u)|,\\
	\Gamma&=\sup_{(x,u)\atop |u|\leq \delta}\sum_{ikm}|\g_{ikm}(x,u)|.
	\ea
	Applying a Gronwall type estimate to the characteristic equation for $w_i$, we get
	\be\label{eq:gronwall1}
	|w_i(x,t)|\leq y(t),
	\ee
	where $y(t)$ satisfies
	\ba\label{eq:gronwall2}
	y'(t)&=Zy(t)+\Gamma y(t)^2+K(d,T)U,\\
	y(0)&=W_0.
	\ea
	Then we claim that there exists a $\overline{W}>0$ and $0<C<\infty$ depending only on $Z$, $\Gamma$, $\delta$, the support of $u(x,0)$ and $t_0$ such that 
	\be\label{eq:gronwall3}
	y(t)\leq CW_0
	\ee 
	holds for all $t\leq t_0$ and all $W_0\leq\overline{W}$. Noting that we have $U\leq \delta$ a priori, we are free to bound $UK(d,T)$ by $\tilde{c}\delta W_0$ for some $\tilde{c}$ sufficiently small. In particular, we can choose $\tilde{c}$ so small that $\tilde{c}\delta W_0\leq ZW_0$. Hence, at time 0 $K(d,T)U\leq Zy(0)$; but since the right hand side of \eqref{eq:gronwall2} is positive, the bound $K(d,T)U\leq Zy(t)$ persists for longer time. This allows us to bound the solution $y(t)$ by the solution $\tilde{y}$ of
	\ba
	\tilde{y}'(t)&=2Z\tilde{y}(t)+\Gamma \tilde{y}(y)^2,\\
	\tilde{y}(0)&=W_0.
	\ea 
	The claim can be proven by writing $Y(t)=\tilde{y}(t)e^{-2Zt}$ for $t\leq t_0$ and noting that $Y(t)$ solves the weighted Riccati equation
	\be
	Y'(t)=(\Gamma e^{2Zt})Y(t)^2.
	\ee 
	Now comparing $Y(t)$ with the Riccati equation with coefficient $\Gamma e^{2Zt_0}$ and following the remainder of the argument for original Riccati equation proves the claim.\\
	
	For $0\leq t\leq t_0$, we recall \eqref{eq:Sest1}
	\be\label{eq:Sest8}
	\frac{d (\b_i(t)-\a_i(t))}{dt}=\l_i(\b_i(t),u(\b_i(t),t))-\l_i(\a_i(t),u(\a_i(t),t)).
	\ee
	We now add and subtract $\l_i(\b_i(t),0)$ and $\l_i(\a_i(t),0)$ to get
	\ba\label{eq:Sest9}
	\Big|\frac{d(\b_i(t)-\a_i(t))}{dt}|\Big|&\leq |\l_i(\b_i(t),u(\b_i(t),t))-\l_i(\b_i(t),0) )|+\\
	&\quad+|\l_i(\b_i(t),0)-\l_i(\a_i(t),0)|+|\l_i(\a_i(t),u(\a_i(t),t))-\l_i(\a_i(t),0)|.
	\ea
	By hypothesis $u(\a_0,0)=u(\b_0,0)=0$, and so we can bound \eqref{eq:Sest9} by
	\ba\label{eq:Sest10}
	\Big|\frac{d(\b_i(t)-\a_i(t))}{dt}|\Big|&\leq C(|u(\b_i(t),t)-u(\b_0,0)|+|u(\a_i(t),t)-u(\a_0,0)|)+\\
	&\quad+|\l_i(\b_i(t),0)-\l_i(\a_i(t),0)|.  
	\ea
	Finally, we are left with controlling $|\l_i(\b_i(t),0)-\l_i(\a_i(t),0)|$. We do this by the fundamental theorem of calculus as follows
	\be\label{eq:Sest11}
	|\l_i(\b_i(t),0)-\l_i(\a_i(t),0)|=\Big|\int_{\a_i(t)}^{\b_i(t)} b_{ii}(s,0)ds\Big|\leq C(e^{-c|\b_i(t)| }-e^{-c|\a_i(t)|} ).
	\ee
	We then see that $|\l_i(\b_i(t),0)-\l_i(\a_i(t),0)|=\cO(1)$ so that for $t\leq t_0$ one has by the characteristic equation for $u$ and the Gronwall type bound on $w_i$ in \eqref{eq:gronwall3}
	\be\label{eq:Sest12}
	\Big|\frac{d(\b_i(t)-\a_i(t))}{dt}\Big|=\cO(1+t_0W_0),
	\ee
	after choosing $d$ large enough that the reaction term in \eqref{eq:ualongichar} satisfies $\delta\cG(d,T)\lesssim W_0$. Hence integrating and using \eqref{eq:Sest7} and \eqref{eq:Sest12} we get
	\be\label{eq:Sest13}
	|\b_i(t)-\a_i(t)|\leq s_0+\cE(d,T)+\cO(t_0+t_0^2W_0)+T\cO(SV+J)=\cO(s_0+\cE(d,T) +TSV+TJ),
	\ee
	as $t_0^2W_0=\cO(s_0^2W_0)=\cO(\theta s_0)$ can be absorbed into $\cO(s_0)$.
\end{proof}
\begin{remark}
	In John's argument, in the final step one assumes that $T\sim \theta^{-1}\gg 1$, which will ensure that all the distance conditions stated at the end of Lemma \ref{lem:USbounds} hold.
\end{remark}
\begin{proposition}\label{prop:Jest}
	We may estimate $J$ by
	\be\label{eq:keyJest}
	J=\cO(s_0W_0+(V+\cE(d,T))TJ+TV(V+\cE(d,T))S+T\cG(d,T)US ),
	\ee
	where we choose the distance $d$ so large that
	\be
	\cE(d,T)\leq\min\{1,s_0,\tilde{c}W_0 \},
	\ee
	holds for $\tilde{c}$ the same small constant as before. We assume that $d$ is so large that
	\be
	\cG(d,T)\leq W_0^2,
	\ee
	holds as well.
\end{proposition}
\begin{proof}
	Following John, we introduce two new quantities
	\ba\label{def:Jauxiliaryquantities}
	\rho_i(z,t)&=\frac{\d X_i}{\d z}(z,t),\\
	v_i(z,t)&=w_i(X_i(z,t),t)\rho_i(z,t).
	\ea
	One should think of $\rho_i$ as some measure of the density of the $i$th characteristics. By a change of variables computation, one discovers that $J$ can be computed in terms of the $v_i$ by
	\be\label{eq:Jcov}
	J=\sup_i \sup_{0\leq t\leq T}\int_{\a_i(t)}^{\b_i(t)}|w_i(x,t)|dx=\sup_i\sup_{0\leq t\leq T}\int_{\a_0}^{\b_0}|v_i(z,t)|dz.
	\ee
	From the characteristic equation for $X_i$ given in \eqref{eq:ithCharacteristicEqn} and the expression for $d\l_i$ in \eqref{eq:KatoFormulasLambdaEta}, we are lead to the evolution equation for $\rho_i$
	\ba\label{eq:rhosystem}
	\frac{\d}{\d t}\rho_i&=\frac{\d}{\d z}\frac{\d X_i}{\d t}=b_{ii}\rho_i+\sum_{m=1}^nc_{iim}w_m\rho_i,\\
	\rho_i(z,0)&=z.
	\ea
	It is clear that $\tilde{\rho}_i(z,t)=0$ for all $t$ is a solution to the first equation of \eqref{eq:rhosystem}. This ensure that our desired solution $\rho_i(z,t)>0$ for all $t$ by uniqueness of solutions, as the initial data is nowhere vanishing. An important consequence of $\rho_i>0$ is that it we may look at the time evolution for $\log(\rho_i)$ given by
	\be\label{eq:logrho1}
	\frac{\d}{\d t}\log(\rho_i)=b_{ii}+\sum_{m=1}^n c_{iim}w_m.
	\ee
	By imposing the additional constraint $\text{dist}([\a_0,\b_0],0)\gtrsim \max\{1,T,|\log(W_0)|\}$, we may ensure that $|b_{ik}(z,t)|\lesssim W_0$ on $[\a_0-\overline{\l}T,\b_0+\overline{\l}T]$ for all $i,k$. From the short time estimate on the $w_i$ and th distance assumption, we find that for $0\leq t\leq t_0$ that
	\be\label{eq:logrho2}
	\frac{\d }{\d t}\log(\rho_i)=\cO(W_0),
	\ee
	which upon integrating in time leads to
	\be\label{eq:logrho3}
	\log(\rho_i)=\cO(t_0W_0)=\cO(\theta).
	\ee
	Exponentiating \eqref{eq:logrho3} leads to the bound for $\rho_i$
	\be\label{eq:shorttimerho}
	\rho_i=\cO(1),
	\ee
	for $0\leq t\leq t_0$.\\
	
	Turning to $v_i$, we first write down the evolution equation as
	\be\label{eq:vevol1}
	\frac{\d v_i}{\d t}=\big(\sum_{k=1}^n \zeta_{ik} w_k+\sum_{k,m}\g_{ikm}w_kw_m+\kappa_i u \big)\rho_i+w_i\big(b_{ii}\rho_i+\sum_{m=1}^nc_{iim}w_m\rho_i \big).
	\ee
	By using the Kronecker $\delta$, \eqref{eq:vevol1} can be more succintly written as
	\be\label{eq:vevol2}
	\frac{\d v_i}{\d t}=\sum_{k=1}^n (\zeta_{ik}+\delta_{ik} b_{ik})\rho_i w_k+\sum_{k,m}(\g_{ikm}+\delta_{ik}c_{ikm} )w_kw_m\rho_i+\kappa_i \rho_i u.
	\ee
	For $0\leq t\leq t_0$, appealing to the Gronwall type estimate \eqref{eq:gronwall3}, the distance condition on the support, and the estimate on $\rho_i$ in \eqref{eq:shorttimerho}, we find that
	\be\label{eq:shorttimev1}
	\frac{\d v_i}{\d t}=\cO(\rho_iW_0^2+\rho_i\cG(d,T)U)=\cO(W_0^2+\cG(d,T)U).
	\ee
	By assumption, $U\leq\delta$ so by choosing the support far enough, $\cG(d,T)U\leq \delta W_0^2$, refining the short time estimate to
	\be\label{eq:shorttimev2}
	\frac{dv_i}{dt}=\cO(W_0^2).
	\ee
	For $t\geq t_0$, we make the crucial observation that $\frac{\d v_i}{\d t}$ does not include $w_i^2$ in the sum over $k,m$ as $\g_{iii}+c_{iii}=0$. Moreover, since $(x,t)=(X_i(z,t),t)$ can only lie in one $\cR_i$ for $t\geq t_0$, we get the estimate on $v_i$ for $t\geq t_0$
	\ba\label{eq:longtimev1}
	\frac{\d v_i}{\d t}(z,t)&=\cO(\cE(d,T)\rho_i(|w_i|+V)+\rho_i|w_i|V+\rho_iV^2+\rho_i \cG(d,T)U)+\\
	&\quad=\cO((V+\cE(d,T))|v_i(z,t)|+(V+\cE(d,T))V\rho_i+\rho_i\cG(d,T)U).
	\ea
	Integrating with respect to time, we can bound $v_i$ by
	\ba\label{eq:longtimev2}
	|v_i(z,t)|&\leq |w_i(z,0)|+C\Big( (V+\cE(d,T))\int_0^T|v_i(z,s)|ds+\\
	&\quad+\big(V(V+\cE(d,T))+\cG(d,T)U\big)\int_0^T\rho_i(z,s)ds+t_0W_0^2\Big),
	\ea
	with the constant $C$ uniform in $z\in[\a_0,\b_0 ]$ and $0\leq t\leq T$. Integrating in the spatial variable $z$, we get after applying Fubini's theorem
	\ba\label{eq:Jest1}
	\int_{\a_0}^{\b_0}|v_i(z,t)|dz&\leq \int_{\a_0}^{\b_0}|w_i(z,0)|dz+C\Big((V+\cE(d,T))\int_0^T\int_{\a_0}^{\b_0}|v_i(z,s)|dzds+\\&\quad+\big(V(V+\cE(d,T))+\cG(d,T)U\big)\int_0^T(\b_i(s)-\a_i(s))ds+s_0^2W_0^2  \Big).
	\ea
	The trivial bound $||f||_{L^1([0,T])}\leq T||f||_{L^\infty([0,T])}$ produces our final estimate for $J$, namely
	\be\label{eq:Jest2}
	J=\cO(s_0W_0+(V+\cE(d,T))TJ+TV(V+\cE(d,T))S+T\cG(d,T)US ),
	\ee
	as $s_0W_0$ small allows us to absorb $s_0^2W_0^2$ into $s_0W_0$.
\end{proof}
The estimate on $J$ that we have looks slightly different from the estimate on $J$ in John's argument, however, we've assumed that $\cE(d,T)\leq W_0$. In the final step, we will start by assuming that $s_0V=\cO(\theta)$, and so $\cE(d,T)$ will at least be a priori comparable to $V$ provided $s_0=\cO(1)$. Note also the presence of the lower order term $T\cG(d,T)US$.\\

The final estimate needed is for the estimate for $V$.
\begin{proposition}\label{prop:Vest}
	Suppose that $\emph{dist}([\a_0,\b_0 ],0)\gtrsim \max\{1,T,|\log(W_0)| \}$ holds. Then one can estimate $V$ by
	\be\label{eq:keyVest}
	V=\cO(s_0W_0^2+(V+\cE(d,T))(VT+(V+\cE(d,T))TJ+TV(V+\cE(d,T))S)+\cG(d,T)TU).
	\ee
	As in Lemma \ref{lem:USbounds} and Proposition \ref{prop:Jest}, we assume that $d\gtrsim\max\{1,T,|\log(W_0)| \}$ so that
	\be
	\cE(d,T)\leq\min\{1,s_0,\tilde{c}W_0,W_0^2\},
	\ee
	and
	\be
	\cG(d,T)\leq W_0^2,
	\ee
	hold.
\end{proposition}
\begin{proof}
	Let $(x,t)\not\in \cR_i$. Then there exists a $z\not\in[\a_0,\b_0]$ such that $(x,t)=(X_i(z,t),t)$. As $z$ is outside the support of $u(x,0)$, we necessarily have $w(z,0)=0$. Moreover, we can assume that $z$ is not too far from $[\a_0,\b_0]$ in the sense that $z\in[\a_0-\overline{\l}T,\b_0+\overline{\l}T ]$, as outside this larger interval $w$ is identically zero. Integrating $w_i$ along the characteristic starting at $z$, we find
	\ba\label{eq:Vest1}
	w_i(x,t)&=\int_0^t\sum_k \zeta_{ik}(X_i(z,\tau),u(X_i(z,\tau),\tau))w_k+\sum_{k,m}\g_{ikm}(X_i(z,\tau),u(X_i(z,\tau),\tau))w_kw_m+\\
	&\quad+\kappa_i(X_i(z,\tau),u(X_i(z,\tau)))u(X_i(z,\tau),\tau) d\tau  
	\ea
	As before, we split into two cases depending on how $t$ compares to $t_0$. If $t\leq t_0$, then we have that the integrand of \eqref{eq:Vest1} is at most $W_0^2$ giving us the short time estimate
	\be\label{eq:shorttimewi}
	|w_i(z,t)|=\cO(t_0W_0^2).
	\ee
	In the other case, we introduce sets $\o_k$ defined by
	\be\label{def:omegak}
	\o_k=\{0\leq\tau\leq T: (X_i(z,\tau),\tau)\in\cR_k \}.
	\ee
	By assumption $\o_i$ is empty. If $\tau\geq t_0$, then $\tau$ is contained in at most one $\cR_k$. For that specific $k$, $|w_k||w_m|\leq V|w_k|$ otherwise $|w_k||w_m|\leq V^2$. For reaction term $\kappa_i u$, one can do little better than
	\be\label{eq:Vreacest}
	\Big|\int_0^t \kappa_i(X_i(z,\tau),u(X_i(z,\tau),\tau))u(X_i(z,\tau),\tau)d\tau\Big|\leq \cG(d,T)TU.
	\ee	For the linear terms, we use the distance on the support to bound them by
	\be\label{eq:linearest1}
	\Big| \int_0^t\sum_{k=1}^n\zeta_{ik}w_kd\tau\Big|\leq \cE(d,T)\sum_k\Big(\int_{\o_k}|w_k|d\tau+\int_{[0,T]\backslash \o_k}|w_k|d\tau \Big). 
	\ee
	Outside of $\o_k$ one has $|w_k|\leq V$, and so \eqref{eq:linearest1} is controlled by
	\be\label{eq:linearest2}
	\Big| \int_0^t\sum_{k=1}^n\zeta_{ik}w_kd\tau\Big|\leq \cE(d,T)VT+\cE(d,T)\sum_{k=1}^n\int_{\o_k}|w_k(X_i(z,\tau),\tau)|d\tau.
	\ee
	We may bound the quadratic terms of \eqref{eq:Vest1} in a similar manner by
	\be\label{eq:quadraticest1}
	\Big|\sum_{k,m}\g_{ikm}(X_i(z,\tau),u(X_i(z,\tau),\tau))w_kw_md\tau \Big|\leq C(V^2T+V\sum_{k=1}^n\int_{\o_k}|w_k(X_i(z,\tau),\tau)|d\tau).
	\ee
	Thus, we have reduced the problem of estimating $V$ to estimating
	\be\label{eq:omegakint1}
	\int_{\o_k}|w_k(X_i(z,\tau),\tau)|d\tau.
	\ee
	We would like to exchange the integral with respect to the time variable $\tau$ for an integral with respect to space, as our setup gives us pointwise control on $w_k$ outside $\cR_k$ and spatial integrals of $w_k$ in $\cR_k$. To exchange the time variable for a suitable space variable in \eqref{eq:omegakint1} fix $k$, and let $\tau\in\omega_k$. Then there exists some $y=y(\tau)$ so that
	\be\label{eq:switchingcharacteristics}
	X_i(z,\tau)=X_k(y(\tau),\tau).
	\ee
	Differentiating \eqref{eq:switchingcharacteristics} with respect to $\tau$, we find that
	\be\label{eq:tauchainrule}
	\l_i=\l_k+\rho_k\frac{d y}{d\tau}.
	\ee
	As $\rho_k>0$ and $\l_k-\l_i$ has a definite sign, $y(\tau)$ is a strictly monotone function of $\tau$ and is thus invertible. Changing variables in \eqref{eq:omegakint1}, we find that
	\be\label{eq:omegakin2}
	\int_{\o_k}|w_k(X_i(z,\tau),\tau)|d\tau=\int_{I_k} \frac{|w_k(X_k(y,\tau),\tau)|\rho_k(y,\tau)}{|\l_i-\l_k| }dy,
	\ee
	for some subinterval $I_k\subset[\a_0,\b_0 ]$. By assumption, $|\l_i-\l_k|$ is bounded from below, so \eqref{eq:omegakin2} can be bounded by
	\be\label{eq:omegakint3}
	\int_{\o_k}|w_k(X_i(z,\tau),\tau)|d\tau\leq C\int_{I_k}|v_k(y,\tau(y))|dy.	
	\ee
	In the proof of the estimate \eqref{eq:keyJest}, we obtained the pointwise bound \eqref{eq:longtimev2} on $v_k$, valid for all $(y,\tau(y))$. Plugging this estimate into \eqref{eq:omegakint3} produces
	\ba\label{eq:omegakint4}
	\int_{\o_k}|w_k(X_i(z,\tau),\tau)|d\tau&\leq C\int_{\a_0}^{\b_0}\Big(|w_i(y,0)|+\Big( (V+\cE(d,T))\int_0^T|v_i(y,s)|ds+\\
	&\quad+V(V+\cE(d,T))\int_0^T\rho_i(y,s)ds+t_0W_0^2\Big) \Big)dy=\\
	&\quad=\cO(s_0W_0+(V+\cE(d,T))TJ+TV(V+\cE(d,T))S). 
	\ea
	We've also expanded $I_k$ into $[\a_0,\b_0]$ as well in \eqref{eq:omegakint4}. At this point, following the same procedure for the $J$ estimate furnishes our final estimate on $V$
	\ba\label{eq:Vest2}
	V&=\underbrace{W_0\cO(VT+s_0W_0+(V+\cE(d,T))TJ+TV(V+\cE(d,T))S)}_{\text{Linear contribution}}+\\
	&\quad+\underbrace{V\cO(VT+s_0W_0+(V+\cE(d,T))TJ+TV(V+\cE(d,T))S)}_{\text{Quadratic contribution} }+\\
	&\quad+\underbrace{\cG(d,T)TU}_{\text{Reaction term}}+\underbrace{\cO(s_0W_0^2)}_{\text{Short time}}= \\
	&\quad=\cO(s_0W_0^2+(V+\cE(d,T))(VT+(V+\cE(d,T))TJ+\\
	&\quad+TV(V+\cE(d,T))S)+\cG(d,T)TU).
	\ea
	Where we've used \eqref{eq:John48} to absorb $s_0W_0V$ into the left hand side of \eqref{eq:Vest2}.
\end{proof}
Putting these results together, we have
\begin{theorem}\label{thm:tableofestimates}
	Suppose $u(x,t)$ is a $C^2$ solution, remaining bounded by $\delta$, to \eqref{eq:perturbedeqn} with initial data $u(x,0)$ compactly supported in $[\a_0,\b_0 ]$ with
	\be\label{eq:distcondition}
	\emph{dist}([\a_0,\b_0],0)\gtrsim\max\{1,T,|\log(s_0)|,|\log(W_0)| \},
	\ee
	and satisfying
	\be
	\max_x s_0^2|u_{xx}(x,0)|=\theta\ll 1,
	\ee
	for $s_0=\b_0-\a_0$.
	Then the following estimates hold
	\begin{subequations}
		\begin{align}
		V&=\cO(s_0W_0^2+(V+\cE(d,T))(VT+(V+\cE(d,T))TJ+\\
		&\quad\nonumber+TV(V+\cE(d,T))S)+\cG(d,T)TU),\\
		U&=\cO(s_0V+TV+J),\\
		S&=\cO(s_0+\cE(d,T)+VTS+JT),\\
		J&=\cO(s_0W_0+(V+\cE(d,T))TJ+TV(V+\cE(d,T))S+T\cG(d,T)US ).
		\end{align}
	\end{subequations}
\end{theorem}
\subsection{Final Steps}
We are now in a position to finish John's argument in showing blowup of some solution. We begin by constructing a suitable initial condition.
\begin{lemma}\label{lem:initialdata}
	For a given index $i$ and all $\theta$ small enough, there exists a nonzero function $f$ which is $C^2$ with $\sup_x|f''(x)|\leq \theta$, compactly supported on an interval of width 1, and such that there exists a $y$ in the support of $f$ so that
	\be
	\eta_i(y,f(y))f'(y)=\sup_i \sup_x |\eta_i(x,f(x))f'(x)|.
	\ee
\end{lemma}
\begin{proof} 
	Let $\phi(x)$ be a scalar function which is smooth, not identically zero, and compactly supported in $[-1/2,1/2]$. For simplicity assume that $||\phi''||_{\infty}=1$. Let $f(x)=f_{i,\theta}(x)=\theta\phi(x-x_0)\xi^i_\infty(0)$, where $\xi^i_\infty(u)$ is given by
	\be
	\xi^i_\infty(u)=\lim_{|x|\to\infty} \xi^i(x,u),
	\ee
	and $x_0\in\RR$ to be determined. By a continuity argument
	\be
	\lim_{|x_0|\to\infty}\eta_i(x_0+t,u)\xi^j_\infty(u)=\delta_{ij},
	\ee
	with uniform convergence on $|t|\leq\frac{1}{2}$ and $|u|\leq\delta$. If we let $\eta_i^\infty(u)$ be defined a manner analogous to $\xi_\infty^i(u)$, then we note that
	\be
	\eta_j^\infty(f(x_0+t))\xi_\infty^i(0)=\delta_{ij}+\cO(\theta),
	\ee
	for all $|t|\leq\frac{1}{2}$. In particular, what we find is that
	\be
	\eta_j(x_0+t,f(x_0+t))\xi^i_\infty(0)\to\eta_j^\infty(\theta \phi(t)\xi_\infty^i(0))\xi_i^\infty(0)=\delta_{ij}+\cO(\theta),
	\ee
	as $|x_0|\to\infty$. In particular, this implies that
	\be
	\eta_j(x,f(x))f'(x)=\theta\phi'(x)\delta_{ij}+\cO(\theta^2),
	\ee
	and so for $\theta\ll 1$, we have that there exists $y$ in $\supp(f)$ such that
	\be
	|\eta_i(y,f(y))f'(y)|=\sup_i \sup_x |\eta_i(x,f(x))f'(x)|,
	\ee
	and by sending $\phi\to-\phi$ if needed, we can remove the absolute value bars on the left hand side.
\end{proof}	
Adapting a result of Lax in \cite{La} as in \cite{B}, one can show that for all small smooth initial data, there is a unique classical solution to \eqref{eq:perturbedeqn} on some time interval with the prescribed initial data. Our solution $u(x,t)$ will then be the unique smooth solution to \eqref{eq:perturbedeqn} with the initial data constructed in Lemma \ref{lem:initialdata}.\\ 

We let $T$ increase from 0 so that the inequalities
\begin{subequations}\label{subeq:originalestimates}
	\begin{align}
	TW_0&\leq \max_i\frac{4}{ \g_{iii}^\infty(0)}=\vartheta,\\
	TV&\leq\sqrt{\theta},\\
	J&\leq\sqrt{\theta},\\
	V&\leq\theta,\\
	U&\leq\sqrt{\theta}.
	\end{align}
\end{subequations}
For $\g_{iii}^\infty(0)$ the value of $\g_{iii}$ with $B=0$, or equivalently,
\be\label{def:gammainfinity}
\g_{iii}^{\infty}(0)=\lim_{|x|\to\infty}\g_{iii}(x,0). 
\ee
Note that these inequalities are all valid at $T=0$ for $\theta$ small enough since there $U=\cO(\theta)$, $T=V=0$, $J=\cO(W_0)=\cO(\theta)$. For the fourth estimate in Theorem \ref{thm:tableofestimates}, we may improve the estimate on $J$ to
\be\label{eq:improvedJest1}
J=\cO(W_0+\sqrt{\theta}J+(T\cE(d,T))J+\sqrt{\theta}(V+\cE(d,T))S+T\cG(d,T)US ).
\ee
As we've chosen $d\sim \theta^{-1}\gg |\log(W_0)|\sim|\log(\theta)|$, we may assume that $\cE(d,T)\leq \min\{1,c\vartheta W_0,s_0\}$ for some constant $c\ll1 $ so small that the term $\cO(\cE(d,T)T)J$ can be absorbed onto the left hand side, that is $\cO(T\cE(d,T))\leq\frac{1}{2}$. As $\theta$ is small, the $\sqrt{\theta}J$ on the right hand side can be absorbed as well. This improves the $J$ estimate to
\be\label{eq:improvedJest2}
J=\cO(W_0+\sqrt{\theta}VS+T\cG(d,T)US).
\ee
To handle the reaction term, we note that $\cG(d,T)\leq cW_0^2$ for some universal constant $c>0$ sufficiently small, hence $T\cG(d,T)U\lesssim W_0U$. This furthers improves the $J$ estimate to
\be\label{eq:improvedJest3}
J=\cO(W_0+\sqrt{\theta}VS+W_0 US)
\ee
From this, we find that $JT$ is of the order
\be\label{eq:improvedJTest}
JT=\cO(1+\theta S+US ).
\ee
Using the bounds $U\leq \sqrt{\theta}$ and \eqref{eq:improvedJTest} allows us to improve the $S$ estimate to
\be\label{eq:improvedSest}
S=\cO(1+\cE(d,T)+VTS+JT)=\cO(1+ \sqrt{\theta}S+\theta S )=\cO(1).
\ee
Plugging \eqref{eq:improvedSest} into \eqref{eq:improvedJest3} produces our final refinement of the $J$ estimate
\be\label{eq:improvedJest4}
J=\cO(W_0+\sqrt{\theta}V+W_0U )=\cO(\theta+\theta^{\frac{3}{2}})=\cO(\theta).
\ee
Next, we use our refined estimates on $S$ and $JT$ in \eqref{eq:improvedSest} and \eqref{eq:improvedJTest} respectively, to refine the estimate on $V$.
\be\label{eq:improvedVest1}
V=\cO(W_0^2+(V+\cE(d,T))(VT+(V+\cE(d,T))+TV(V+\cE(d,T)))+\cG(d,T)TU).
\ee
Applying the bound $TV\leq\sqrt{\theta}$ and using $\cO(T\cE(d,T)))\leq\frac{1}{2}$ as in the $J$ estimate, we get a further refinement of the $V$ estimate
\be\label{eq:improvedVest2}
V=\cO(W_0^2+\sqrt{\theta}V+(V+\cE(d,T))^2+\sqrt{\theta}(V+\cE(d,T))^2+\cG(d,T)TU)+\frac{1}{2}V.
\ee
This then gives
\be\label{eq:improvedVest3}
V=\cO(W_0^2+(V+\cE(d,T))^2+\cG(d,T)TU).
\ee
A priori, we have $V,W_0=\cO(\theta)$ and $\cE(d,T)\lesssim W_0$, allowing us to absorb the $\cE(d,T)V$ and $V^2$ terms of \eqref{eq:improvedVest3} into the left hand side; hence we conclude
\be\label{eq:improvedVest4}
V=\cO(W_0^2+\cG(d,T)TU)=\cO(W_0^2+W_0U).
\ee
From \eqref{eq:improvedVest3}, we may improve the estimate on $TV$ to
\be\label{eq:improvedTVest1}
TV=\cO(TW_0^2+\cG(d,T)T^2U)=\cO(TW_0^2+T^2W_0^2U).
\ee
Using $TW_0\leq\vartheta=\cO(1)$, we get the refinement
\be\label{eq:improvedTVest2}
TV=\cO(TW_0^2+T^2W_0^2U)=\cO(W_0+U).
\ee
The final estimate to improve is the $U$ estimate. The improvement is obtained by using \eqref{eq:improvedVest4}, \eqref{eq:improvedTVest2}, and \eqref{eq:improvedJest3}
\be\label{eq:improvedUest}
U=\cO(s_0V+TV+J)=\cO(W_0^2+W_0U+W_0+cU+J )=\cO(\theta),
\ee
by using the small constant in $\cG(d,T)\leq cW_0^2$ to move the $U$ on the right hand side to the left hand side of \eqref{eq:improvedUest}.
To summarize, we've bounced the estimates in \eqref{subeq:originalestimates} off of each other in order to get the better estimates
\begin{subequations}\label{subeq:improvedestimates}
	\begin{align}
	TV&\leq C\theta,\\
	J&\leq C\theta,\\
	V&\leq C\theta^2,\\
	U&\leq C\theta.
	\end{align}
\end{subequations}
for some constant $C$ depending only on the matrices $A$, $B$, the parameters $\delta$ and $\e$ for $\theta$ small enough under the assumptions that $d$ is sufficiently large and $s_0=1$.\\

Suppose $u(x,t)$ remains $C^2$ for all $0\leq t\leq T$ for $T=\vartheta W_0^{-1}$ and that for some $i$ there exists a $z\in[x_0-\frac{1}{2},x_0+\frac{1}{2}]$ such that
\be\label{eq:blowupassumption}
w_i(z,0)=W_0.
\ee
We consider the evolution of $w_i(z,0)$ along the characteristic $X_i(z,t)$, that is we look at the function $w(t)$ defined by
\be
w(t)=w_i(X_i(z,t),t).
\ee
Appealing to the characteristic equation for $w(t)$, the estimate \eqref{eq:gronwall3}, and the assumption $|x_0|\gtrsim |\log(W_0)|$ we find that
\be\label{eq:shortevolution1}
\Big|\frac{dw}{dt}\Big|\leq C(\cE(d,T)W_0+W_0^2)\leq CW_0^2,
\ee
for all $0\leq t\leq t_0$. Note that the reaction term $\kappa_i u$ has been absorbed into the $W_0^2$ by using $U\cG(d,T)\leq W_0^2$. Integrating with respect to time leads to the estimate
\be\label{eq:shortevolution2}
|w(t)-W_0|=|w(t)-w(0)|=\cO(t_0 W_0^2)=\cO(W_0^2)=\cO(\theta W_0),
\ee
where we've used $t_0=\cO(s_0)=\cO(1)$ and $W_0=\cO(\theta)$. In particular, provided that $\theta$ is small enough, $w(t)$ will satisfy
\be\label{eq:shortevolution3}
w(t)>\frac{3}{4}W_0 \quad\quad \text{for} \quad\quad 0\leq t\leq t_0.
\ee
For each $0\leq t\leq T$ and all $x$, we have by a mean value theorem type estimate that
\be\label{eq:MVT1}
|\g_{iii}(x,u(x,t))-\g_{iii}(x,0)|=\cO(U)=\cO(\theta).
\ee
In addition, since the perturbation $B(x,u)$ decays exponentially in $x$, we see that by the fundamental theorem of calculus
\be\label{eq:MVT2}
\g_{iii}(x,0)=\g_{iii}^\infty(0)+\int_{-\infty}^x \frac{\d \g_{iii}(y,0)}{\d y}dy=\g_{iii}^\infty(0)+\cE(d,T).
\ee
Here, we've implicitly assumed that $[\a_0,\b_0]\subset (-\infty,0)$ for the sake of definiteness, but an entirely analogous calculation works in the other case and will thus be omitted. In either case, we have that
\be\label{eq:MVT3}
|\g_{iii}(x,u(x,t))-\g_{iii}^\infty(0)|=\cO(U+\cE(d,T))=\cO(\theta).
\ee
This implies that, for $\theta$ small enough and all $x$ in the support of $u$ that
\be
\g_{iii}(x,u(x,t))>\frac{1}{2} \g_{iii}^\infty(0).
\ee
For $k\not=i$ and $t>t_0$ we have that $|w_k(X_i(z,t),t)|\leq V$. Using this information and \eqref{eq:MVT3} in the characteristic equation for $w$, we have
\be\label{eq:wcharest1}
\frac{dw}{dt}>\frac{1}{2}\g_{iii}^\infty(0)w^2-\cE(d,T)|w|-n\cE(d,T)V-\Gamma(V|w|+V^2)-K(d,T)U,
\ee
for $\Gamma$ and $K(d,T)$ as in \eqref{def:gronwallparam}. At $t=t_0$, we can use $\cE(d,T)\ll W_0\lesssim w(t_0)$, $K(d,T)\lesssim W_0^2$, and $V=\cO(\theta W_0)$ to conclude that for $\theta$ small enough, $\frac{d w}{dt}(t_0)>0$. Hence $w(t)$ is increasing which implies that $V=\cO(\theta |w|)$ persists for larger times. For $\theta$ small enough, we then have for $t>t_0$
\be\label{eq:wcharest2}
\frac{dw}{dt}>\frac{3}{8}\g_{iii}^\infty(0)w^2.
\ee
However, the solution to this Riccati equations blows up at some time $t$ satisfying
\be
t\leq t_0+\frac{8}{3\g_{iii}^\infty(0)w(t_0)}\leq t_0+\frac{2}{3}T<T.
\ee
This final conclusion comes from the observation that
\be
t_0=\cO(\frac{\theta}{W_0})=\cO(\theta T)<\frac{1}{3}T,
\ee
for $\theta$ sufficiently small. As $|w_i(x,t)|=|\eta_i(x,u(x,t))w(x,t)|\leq |w(x,t)|$ by Cauchy-Schwarz, we see that the $L^\infty$ norm of $w$ blows up in finite time. This concludes the adaptation of John's argument incorporating the effects of exponentially small perturbations. We make some final concluding remarks. The first is that exponential decay in the perturbation is not essential here, so long as the perturbation $B(x,u)$ is smooth sufficiently far away from $x=0$ and is $L^1$ with $L^1$ derivative, then the argument can be adapted at the price of losing the ability to estimate how far the support of $u(x,0)$ needs to be from zero. Another remark is that there only needs to be one genuinely nonlinear field in this argument provided that the initial data is largest in the genuinely nonlinear direction, in the sense that the value $W_0$ is achieved by some $w_i$ for which the corresponding coefficient $\g_{iii}\not=0$. This is because we only showed that one specific $w_i$ blows up in finite time, we have little control over what the other $w_k$ are doing. 
H\"ormander in \cite{H} does a more refined analysis of ODE of the form
\be
\frac{dw}{dt}=a_0(t)w^2+a_1(t)w+a_2(t),		
\ee
to show that the blowup time is asymptotically determined by the corresponding blowup time for the Riccati equation
\be
\frac{dw}{dt}=\g_{iii}^\infty(0)w^2.
\ee 
As our characteristic equation is of the same form as the the ODE \cite{H} studies, one can adapt the methods to show the same result here by showing the same bounds as in \cite{H} hold for the characteristic equation obtained here. An alternative argument is to note that $V=\cO(\theta W_0)$, $\cE(d,T)=\cO(W_0)$, $\cG(d,T)=\cO(W_0^2)$, and $U=\cO(\theta)$ can be combined to show that
\be
(\g_{iii}(x,u(x,t))-\cO(\theta))w^2\leq\frac{dw}{dt}\leq(\g_{iii}(x,u(x,t))+\cO(\theta))w^2,
\ee
for $t\geq t_0$. The final step is to note that $\g_{iii}(x,u(x,t))\to \g_{iii}^\infty(0)$ locally uniformly as $\theta\to0$ by $U=\cO(\theta)$ and the distance condition on the support of $u$.\\

The final remark we make here is that John's argument shows that a large class of initial data leads to blow up in the solution, our adaption does not lead to as general a result. It would be interesting to see if generic small data supported far from $x=0$ leads to blow up in this perturbed setup.

\subsection{Application to the ZND model}
Recall the ZND model
\ba
	v_t-u_x&=0,\\
	u_t+p(v,E)_x&=0,\\
	E_t+(pu)_x&=qk\phi(T)z,\\
	z_t&=-k\phi(T)z,
\ea
where $\phi(T)$ is 1 for $T\geq T_i$ and 0 for $T<T_i$. Consider a shock solution $\ul{U}=(\ul{v},\ul{u},\ul{E},\ul{z})$ with shock speed $\s\not=0$ and such that the temperature $T(\ul{U}(x))$ satisfies
\be
	\inf_{x<0}T(\ul{U}(x))>T_i.
\ee
Although the ZND model is not of the form \eqref{eq:unperturbedeqn}, we will still be able to show blowup using the prior method for suitable perturbations of the shock $\ul{U}$.
Making a Galilean change of coordinates into the frame where the shock $\ul{U}$ is stationary, we find that $U=(v,u,E,z)$ satisfies
\ba
	v_t-\s v_x-u_x&=0,\\
	u_t-\s u_x+p_x&=0,\\
	E_t-\s E_x+(pu)_x&=q k\phi(T)z,\\
	z_t-\s z_x&=-k\phi(T)z.
\ea
Now consider a solution of the ZND model of the form $U=\ul{U}+\hat{U}$, where the support of $\hat{U}(x,0)=[\a_0,\b_0]$ is a compact subset of $(-\infty,0)$ with $|\b_0|$ satisfying the distance conditions in Theorem \ref{thm:tableofestimates}. Write $\hat{U}=(\hat{v},\hat{u},\hat{E},\hat{z})$ and further assume that $\hat{z}(x,0)$ is identically 0. We will show that this assumption allows us to eliminate $z$ from the ZND model, reducing the system to gas dynamics. The importance of this reduction is that gas dynamics is a system of the form \eqref{eq:unperturbedeqn}. If the solution $\hat{U}$ remains sufficiently small in $L^\infty$ so that the temperature remains above $T_i$ for all time, then $\hat{z}$ will satisfy
\be\label{eq:hatzeqn}
	\hat{z}_t-\s\hat{z}_x=-k\hat{z},
\ee
as the distance condition on $[\a_0,\b_0]$ ensures that $\supp(\hat{U}(x,t))\subset(-\infty,0)$ for all $t$ up to $T_*$. Solving \eqref{eq:hatzeqn} with the initial data $\hat{z}(x,0)=0$ shows that $\hat{z}(x,t)=0$ for all time. Looking at the $S$ equation in the ZND model, and writing $E=\ul{E}+\hat{E}$ and subtracting off the equation for $\ul{E}$ we find
\be\label{eq:hatEeqn}
	\hat{E}_t-\s\hat{E}_x+(pu-\ul{pu})_x=0.
\ee
We have eliminated $z$ from the ZND model for perturbations $\hat{U}$ of the form $\hat{U}=(\hat{v},\hat{u},\hat{E},0)$, leaving us with a system of conservation laws known to have at least one genuinely nonlinear field.
\begin{remark}
	Generic shocks in the ZND model do not exponentially converge to zero as $x\to-\infty$ in every component, but that is easily remedied by writing $\ul{U}=\ul{U}(-\infty)+\tilde{\ul{U}}$, adjusting the pressure function $p$ appropriately, with $\tilde{\ul{U}}$ and all of its derivatives decaying exponentially as $x\to-\infty$.\\
	
	There are shock solutions in the sonic case, called Chapman-Jouget waves, of the ZND model which do not decay exponentially in $x$; but instead exhibit power law decay with $\ul{U}-\ul{U}(-\infty)$ decaying like $x^{-1}$.
\end{remark}
We also note that this argument works equally for the Majda model, showing the necessity of the weight in the stability result Theorem \ref{thm:MajdaStab}.
\subsection{Proof of the corollary}
\begin{corollary}
For every $C>0$, $\theta >0$, $\delta >0$, $s>\dfrac{3}{2}$ and $\epsilon >0$ there exists $\phi$ smooth and compactly supported away from the shock such that $\|\phi\|_{L^2} \leq \delta$ and the solution obtained is given with a perturbation $\Phi$ from the initial wave it gives raise to is defined on some time interval $[0,T]$ and it satisfies $\|U(t)\|_{L^2} \leq \dfrac{C\epsilon}{\theta}$ on $[0,T]$, and $\|U(T)\|_{H^s} > \epsilon$. \\

In particular, there can not be any damping estimates of the form: there exists $C,\theta,\delta,\epsilon$ all positive and $s > \dfrac{3}{2}$ such that, for every solution $\Phi$ that can be written as $\Phi(t,\cdot)=U(\cdot-\psi(t))+v(t,\cdot)$ with $\|v(t)\|_{H^s} \leq \epsilon$ on $[0,T_0]$ then, on $[0,T_0]$ $$\|v(t)\|_{H^s} \leq Ce^{-\theta t}\|v(0)\|_{H^s} + \int_0^tCe^{-\theta (t-s)}\|v(s)\|_{L^2}ds.$$
Furthermore, there is no orbital stability in $H^s(\R^*)$.
\end{corollary}
\begin{proof}
We constructed before positive constants $C'>0$, $\delta$ and $\mu>0$ and a sequence of smooth compactly supported away from the shock initial perturbation $(v_n)_n$ such that for every $s > \dfrac{3}{2}$ $\|v_n\|_{H^s(\R^*)} = \mathcal{O} \left( \dfrac{1}{n+1} \right) $ and the same is true for the $L^{\infty}$ norm of $v_n$ and $(v_n)_x$, with a support which Lebesgue measure is at most $1$ and the associated solutions (the sequence $(U_n)_n$) to the original equation, with initial data the profile of the wave perturbed by $v_n$ can be written as $U_n(t,\cdot-\sigma t)=\underline{U}(\cdot) + w_n(t,\cdot)$, $w_n$ being in $H^s(\R^*)$ satisfying, on the interval of existence $\|w_n\|_{L^{\infty}} \leq C' \|v_n\|_{L^{\infty}}$, $\supp(w_n(t)) \subset \supp(v_n)+[-C't,C't]$ and $\|(w_n)_x(t)\|_{L^{\infty}} \geq \dfrac{\delta}{n\mu-t}$. \\

Thus, given $\epsilon >0$, for some $t_n \leq n\mu$, we have $\|w_n(t_n)\|_{H^s} > \epsilon$ as a lower bound, and, also, the $L^2$ norm of $w_n(t)$ (for $0\leq t \leq t_n$) is bounded by above by $\|w_n(t)\|_{L^{\infty}} \lambda_1(\supp(w_n(t)))$, and thus $\|w_n(t)\|_{L^2} \leq \dfrac{C\sqrt{C}(\sqrt{t_n}+\sqrt{\lambda_1(\supp(v_n))})}{n+1}$, thus, up to some multiplicative constant that does not depend on $n$, it is of size at most $\sqrt{n+1}^{-1}$. \\

Hence, given $C>0$ and $\theta >0$, for $n$ big enough, we have, $\|v_n\|_{H^s} \leq \dfrac{\epsilon}{2}$, and, for $t_n$ as before, that on $[0,t_n]$, $U_n$ remains a Lax solution, and $\|w_n\|_{L^2} \leq \dfrac{\epsilon \theta}{2C}$ and $\|w_n(t_n)\|_{H^s} > \epsilon$ thus making the damping estimate described before impossible with these constants $C$, $\epsilon$ and $\theta$. Any damping estimate of this form is, thus, impossible. \\

Furthermore, $(v_n)_n$ goes to $0$ in $H^s$ (for every fixed $s$) as $n$ goes to $+\infty$, but$(U_n)_n$ is not even globally defined. Thus, it precludes orbital stability results in $H^s$.
\end{proof}
\subsection{The case of weighted norms}
In the following, we prove an instability result for weighted norms as used in the Majda model.
\begin{lemma}
We consider the framework considered before of \begin{align*}
u_t+(A(u)+B(x,u))u_x+G(x,u)u=0
\end{align*}
satisfying the same assumptions, and the added constraint that, furthermore, $\Lambda_n(0)<0$. Let $\alpha >0$. There exists $\epsilon_0 >0$ such that, for every $\epsilon \in (0,\epsilon_0)$, for all $\delta >0$ there exists some initial data $u_0 \in C^{\infty}_c(\R)$ such that $\|u_0\|_{H^s_{\alpha}}\leq \delta$ and the solution $u$ has its $H^2_{\alpha}$ norm that gets bigger than $\epsilon$ on the interval of existence of the solution.
\end{lemma}
\begin{proof}
Under the assumptions of the lemma, with $\alpha > 0$ and $s > \dfrac{3}{2}$, let $\xi_n$ be a right eigenvector of $A(0)$ associated with the eigenvalue $\Lambda_n(0)$. We also recall that $\|B(x,u)u_x\|_{L^2_{\alpha}} \lesssim \|u\|_{L^{\infty}}\|u_x\|_{L^2_{\alpha}}$ and $\|G(x,u)u\|_{L^2_{\alpha}} \lesssim (\sup_{x \leq x_0}\|\underline{u}_x\|)\|u\|_{L^2_{\alpha}}$ if $\supp(u) \subset (-\infty,x_0]$. We define for a given solution $u$, $N(u) : (t,x) \mapsto -G(x,u)u - B(x,u)u_x$. With $\phi$ a smooth function with support in $[-1,0]$ that is not the zero function, we have that $v_p : x \mapsto c_p e^{-\alpha p}\phi(\cdot+p) \xi_n$ goes to $0$ in $H^2_{\alpha}$ with $c_p$ a sequence of positive numbers going to $0$ that will be specified later, and, furthermore, it induces a solution $u_p$ defined on $[0,T)$ such that $w_p:=\eta_n \cdot u_p$ satisfies $(w_p)_t + \lambda_n(w_p)_x = \eta_n \cdot (A(0)-A(u)-B(x,u))(u_p)_x - \eta_n \cdot G(x,u_p) u_p$. \\

Assume, by contradiction, that for all $\epsilon$, there exists $p_\epsilon \in \N$ such that, for every $p \geq p_\epsilon$ $u_p$ stays in the ball of radius $\epsilon$ and center $0$ in  $H^1$ for every $t \geq 0$ (and, thus, that the solution is defined on $\R^+$). In particular, given $T \geq 0$, we thus have $$w_p(t,x)=w_p(t,x-\lambda_n(t-T))+ \int_T^t\eta_n\cdot N(u_p)(s,x-\lambda_n (t-s))ds$$ for every $t \geq T$. \\

Let $$\tau : (-\infty,0] \rightarrow \R_+,$$ such that for $x \leq 0$ we have $\tau(x)=\inf_{y \leq x} \| \underline{u}_x \|$ (with the norm being the usual sup norm on $\R^n$). \\

We want to ensure that the Duhamel term stays smaller than $\dfrac{e^{\alpha T}c_p\|\phi\|_{L^2_{\alpha}}}{2}$ with $T$ such that $\dfrac{e^{\alpha T}c_p\|\phi\|_{L^2_{\alpha}}}{2} > \epsilon$. We will use that $\supp(w_p(t,\cdot)) \subset (-\infty,1-p+\mu_1t]$. \\

Thus, as $T_p=2\dfrac{\ln(\epsilon)-\ln(c_p)+\ln(2)-\ln \left( \|\phi\|_{L^2_{\alpha}} \right) }{\alpha}$ is big enough for our needs, we only need to choose a sequence $(c_p)_p$ that converges slowly enough to $0$ to ensure that $\tau(-p + 1 + \mu_1 T_p)$ goes to $0$ fast enough in comparison. Assuming that $c_p=o(p^{-1})$, we have that $T_p=o(p)$, and thus $\tau(-p + 1 + \mu_1 T_p)\leq \tau(-0.5p)$ for large $p$. As a consequence, we obtain the bounds \begin{align*}
|\eta_n \cdot ((A(u_p)-A(0))(u_p)_x)| \lesssim |u_p||(u_p)_x|
\end{align*}
\begin{align*}
|\eta_n \cdot (B(x,u_p)(u_p)_x)| \lesssim \tau(-0.5p) |(u_p)_x|
\end{align*}
\begin{align*}
|\eta_n \cdot (G(x,u_p)u_p)| \lesssim \tau(-0.5p)|u_p|
\end{align*}
Thus, as long as $\tau(-0.5p)=o(c_p)$, for example $c_p=\max(p^{-1},\tau(-0.5p))^2$, we obtain that there is a contradiction for $p$ big enough, as $\|u_p(T_p,\cdot)\|_{L^{\infty}} \lesssim \epsilon e^{-\dfrac{\alpha p}{2}}$. \begin{align*}
|w_p(T_p,x)| \geq \phi(x-\lambda_nT_p+p)c_pe^{-\alpha p} - \int_0^{T_p}|N(t,x-\lambda_n(T_p-t))|dt
\end{align*}
Hence $$\|w_p(T_p,\cdot)\|_{L_{\alpha}^2} \geq \|\phi\|_{L^2_{\alpha}}e^{(p-\lambda_nT_p)\alpha}e^{-p\alpha }-\int_0^{T_p}\|N(t,x-(T_p-t)\lambda_n)\|_{L^2_{\alpha}}dt$$
And so $$\|w_p(T_p,\cdot)\|_{L_{\alpha}^2} \geq \dfrac{\|\phi\|_{L^2_{\alpha}}}{2}e^{-\lambda_nT_p\alpha}$$
for $p$ big enough.
\end{proof}
We now need to check that the lemma applies to ZND. We will first study the eigenvalues of $A(V)$ for a given $V$, for the reduced problem obtained in subsection $3.5$.
$$\begin{pmatrix}
0 & -1 & 0 \\
p_v & p_u & p_E \\
u p_v & u p_u + p & u p_E
\end{pmatrix}$$
which, as its kernel is nontrivial and under the assumption that $v>0$, $u>0$, $E>0$, $(p_v(v,u,E),p_E(v,u,E))\neq (0,0)$ and $p(v,u,E)>0$ we have that the matrix is diagonalizable with real eigenvalues if and only if $pp_E>p_v$ in which case the spectrum of the matrix is $\left\lbrace -\sqrt{pp_E-p_v},0,\sqrt{pp_E-p_v} \right\rbrace$. As $\sigma$ is positive, we have $\sigma > \Lambda_3(V)$, which allow us to apply the previous lemma to this problem, and we will do it under the assumption that the solution studied is close to a shock as described before and that we assume to be admissible, in the sense that it satisfies the assumptions of Theorem \ref{th1}.


\bibliographystyle{alphaabbr}
\bibliography{biblio_chocs_2}

\newcommand{\etalchar}[1]{$^{#1}$}
\newcommand{\SortNoop}[1]{}
\begin{thebibliography}{JNR{\etalchar{+}}19}

\bibitem[Ali95]{Al}
S.~Alinhac.
\newblock {\em Blowup for nonlinear hyperbolic equations}, volume~17 of {\em
  Progress in Nonlinear Differential Equations and their Applications}.
\newblock Birkh\"{a}user Boston, Inc., Boston, aj2, 1995.

\bibitem[B{\"a}r22]{B}
J.~B{\"a}rlin.
\newblock Formation of singularities in solutions to nonlinear hyperbolic
  systems with general source terms.
\newblock {\em arXiv preprint arXiv:2204.12101}, 2022.

\bibitem[BGS07]{BS}
S.~Benzoni-Gavage and D.~Serre.
\newblock {\em Multidimensional hyperbolic partial differential equations}.
\newblock Oxford Mathematical Monographs. The Clarendon Press, Oxford
  University Press, Oxford, 2007.
\newblock First-order systems and applications.

\bibitem[BR22]{BR}
P.~Blochas and L.~M. Rodrigues.
\newblock Uniform asymptotic stability for convection-reaction-diffusion.
\newblock {\em arXiv preprint arXiv:2201.13436}, 2022.

\bibitem[Bre00]{Bress}
A.~Bressan.
\newblock {\em Hyperbolic systems of conservation laws}, volume~20 of {\em
  Oxford Lecture Series in Mathematics and its Applications}.
\newblock Oxford University Press, Oxford, 2000.
\newblock The one-dimensional Cauchy problem.

\bibitem[BCP00]{BCP}
A.~Bressan, G.~Crasta, and B.~Piccoli.
\newblock Well-posedness of the {C}auchy problem for {$n\times n$} systems of
  conservation laws.
\newblock {\em Mem. Amer. Math. Soc.}, 146(694):viii+134, 2000.

\bibitem[BLY99]{BLY}
A.~Bressan, T.-P. Liu, and T.~Yang.
\newblock {$L^1$} stability estimates for {$n\times n$} conservation laws.
\newblock {\em Arch. Ration. Mech. Anal.}, 149(1):1--22, 1999.

\bibitem[CG22]{CG}
S.~Chaturvedi and C.~Graham.
\newblock The inviscid limit of viscous burgers at nondegenerate shock
  formation.
\newblock {\em arXiv preprint arXiv:2204.01170}, 2022.

\bibitem[Chr07]{C}
D.~Christodoulou.
\newblock {\em The formation of shocks in 3-dimensional fluids}.
\newblock EMS Monographs in Mathematics. European Mathematical Society (EMS),
  Z\"{u}rich, 2007.

\bibitem[DF79]{FD}
W.~C. Davis and W.~Fickett.
\newblock Detonation: Theory and experiment.
\newblock 1979.

\bibitem[DR20]{DR1}
V.~Duch{\^ e}ne and L.~M. Rodrigues.
\newblock Large-time asymptotic stability of {R}iemann shocks of scalar balance
  laws.
\newblock {\em SIAM J. Math. Anal.}, 52(1):792--820 889, 2020.

\bibitem[DRar]{DR2}
V.~Duch{\^ e}ne and L.~M. Rodrigues.
\newblock Stability and instability in scalar balance laws: fronts and periodic
  waves.
\newblock {\em Anal. PDE}, to appear.

\bibitem[Erp]{Er1}
J.~J. Erpenbeck.
\newblock Stability of steady-state equilibrium detonations.
\newblock {\em Physics of Fluids (U.S.)}.

\bibitem[FR22]{FR}
G.~Faye and L.~M. Rodrigues.
\newblock Exponential asymptotic stability of riemann shocks of hyperbolic
  systems of balance laws.
\newblock {\em arXiv preprint arXiv:2207.12686}, 2022.

\bibitem[Fic79]{F}
W.~Fickett.
\newblock Detonation in miniature.
\newblock {\em American Journal of Physics}, 47:1050--1059, 1979.

\bibitem[FS04]{FS}
S.~Friedlander and D.~Serre, editors.
\newblock {\em Handbook of mathematical fluid dynamics. {V}ol. {III}}.
\newblock North-Holland, Amsterdam, 2004.

\bibitem[GS93]{GS}
I.~Gasser and P.~Szmolyan.
\newblock A geometric singular perturbation analysis of detonation and
  deflagration waves.
\newblock {\em SIAM J. Math. Anal.}, 24(4):968--986, 1993.

\bibitem[GX92]{GX}
J.~Goodman and Z.~P. Xin.
\newblock Viscous limits for piecewise smooth solutions to systems of
  conservation laws.
\newblock {\em Arch. Rational Mech. Anal.}, 121(3):235--265, 1992.

\bibitem[H{\"o}r87]{H}
L.~H{\"o}rmander.
\newblock The lifespan of classical solutions of nonlinear hyperbolic
  equations.
\newblock In {\em Pseudodifferential operators ({O}berwolfach, 1986)}, volume
  1256 of {\em Lecture Notes in Math.}, pages 214--280. Springer, Berlin, 1987.

\bibitem[JLW05]{JLW}
H.~K. Jenssen, G.~Lyng, and M.~Williams.
\newblock Equivalence of low-frequency stability conditions for
  multidimensional detonations in three models of combustion.
\newblock {\em Indiana Univ. Math. J.}, 54(1):1--64, 2005.

\bibitem[Joh74]{J}
F.~John.
\newblock Formation of singularities in one-dimensional nonlinear wave
  propagation.
\newblock {\em Comm. Pure Appl. Math.}, 27:377--405, 1974.

\bibitem[JNR{\etalchar{+}}19]{JNRYZ}
M.~A. Johnson, P.~Noble, L.~M. Rodrigues, Z.~Yang, and K.~Zumbrun.
\newblock Spectral stability of inviscid roll waves.
\newblock {\em Comm. Math. Phys.}, 367(1):265--316, 2019.

\bibitem[JNRZ14]{JNRZ}
M.~A. Johnson, P.~Noble, L.~M. Rodrigues, and K.~Zumbrun.
\newblock Behavior of periodic solutions of viscous conservation laws under
  localized and nonlocalized perturbations.
\newblock {\em Invent. Math.}, 197(1):115--213, 2014.

\bibitem[JZN11]{JZN}
M.~A. Johnson, K.~Zumbrun, and P.~Noble.
\newblock Nonlinear stability of viscous roll waves.
\newblock {\em SIAM J. Math. Anal.}, 43(2):577--611, 2011.

\bibitem[JYZ21]{JYanZ}
S.~Jung, Z.~Yang, and K.~Zumbrun.
\newblock Stability of strong detonation waves for {M}ajda's model with general
  ignition functions.
\newblock {\em Quart. Appl. Math.}, 79(2):357--365, 2021.

\bibitem[JY12]{JYao}
S.~Jung and J.~Yao.
\newblock Stability of {ZND} detonations for {M}ajda's model.
\newblock {\em Quart. Appl. Math.}, 70(1):69--76, 2012.

\bibitem[Kat76]{K}
T.~Kato.
\newblock {\em Perturbation theory for linear operators}.
\newblock Springer-Verlag, Berlin, second edition, 1976.
\newblock Grundlehren der Mathematischen Wissenschaften, Band 132.

\bibitem[Kaw84]{Kaw}
S.~Kawashima.
\newblock Systems of a hyperbolic-parabolic composite type, with applications
  to the equations of magnetohydrodynamics.
\newblock 1984.

\bibitem[KS88]{KSh}
S.~Kawashima and Y.~Shizuta.
\newblock On the normal form of the symmetric hyperbolic-parabolic systems
  associated with the conservation laws.
\newblock {\em Tohoku Math. J. (2)}, 40(3):449--464, 1988.

\bibitem[Lai19]{Lai}
G.~Lai.
\newblock Detonation wave solution to a 1{D} piston problem for the
  {Z}eldovich-von {N}eumann-{D}\"{o}ring combustion model.
\newblock {\em J. Differential Equations}, 267(9):4949--4974, 2019.

\bibitem[Lax73]{La}
P.~D. Lax.
\newblock {\em Hyperbolic systems of conservation laws and the mathematical
  theory of shock waves}.
\newblock Conference Board of the Mathematical Sciences Regional Conference
  Series in Applied Mathematics, No. 11. Society for Industrial and Applied
  Mathematics, Philadelphia, Pa., 1973.

\bibitem[Lev92]{Le}
A.~Levy.
\newblock On {M}ajda's model for dynamic combustion.
\newblock {\em Comm. Partial Differential Equations}, 17(3-4):657--698, 1992.

\bibitem[LXY22]{LXY}
J.~Li, G.~Xu, and H.~Yin.
\newblock On the blowup mechanism of smooth solutions to 1d quasilinear
  strictly hyperbolic systems with large initial data, 2022.

\bibitem[Liu77]{Li3}
T.~P. Liu.
\newblock Decay to {$N$}-waves of solutions of general systems of nonlinear
  hyperbolic conservation laws.
\newblock {\em Comm. Pure Appl. Math.}, 30(5):586--611, 1977.

\bibitem[Liu79]{Li1}
T.~P. Liu.
\newblock Development of singularities in the nonlinear waves for quasilinear
  hyperbolic partial differential equations.
\newblock {\em J. Differential Equations}, 33(1):92--111, 1979.

\bibitem[Liu87]{Li2}
T.-P. Liu.
\newblock Hyperbolic conservation laws with relaxation.
\newblock {\em Comm. Math. Phys.}, 108(1):153--175, 1987.

\bibitem[Liu20]{Li4}
T.-P. Liu.
\newblock {$N$}-waves for conservation laws with linear defeneracy.
\newblock {\em Bull. Inst. Math. Acad. Sin. (N.S.)}, 15(3):187--216, 2020.

\bibitem[LY99]{LiY}
T.-P. Liu and S.-H. Yu.
\newblock Nonlinear stability of weak detonation waves for a combustion model.
\newblock {\em Comm. Math. Phys.}, 204(3):551--586, 1999.

\bibitem[LRTZ07]{LyRaTZ}
G.~Lyng, M.~Raoofi, B.~Texier, and K.~Zumbrun.
\newblock Pointwise {G}reen function bounds and stability of combustion waves.
\newblock {\em J. Differential Equations}, 233(2):654--698, 2007.

\bibitem[LZ04]{LyZ}
G.~Lyng and K.~Zumbrun.
\newblock One-dimensional stability of viscous strong detonation waves.
\newblock {\em Arch. Ration. Mech. Anal.}, 173(2):213--277, 2004.

\bibitem[Maj81]{M}
A.~Majda.
\newblock A qualitative model for dynamic combustion.
\newblock {\em SIAM J. Appl. Math.}, 41(1):70--93, 1981.

\bibitem[Maj83a]{Maj2}
A.~Majda.
\newblock The existence of multidimensional shock fronts.
\newblock {\em Mem. Amer. Math. Soc.}, 43(281):v+93, 1983.

\bibitem[Maj83b]{Maj1}
A.~Majda.
\newblock The stability of multidimensional shock fronts.
\newblock {\em Mem. Amer. Math. Soc.}, 41(275):iv+95, 1983.

\bibitem[MZ05]{MZ}
C.~Mascia and K.~Zumbrun.
\newblock Stability of large-amplitude shock profiles of general relaxation
  systems.
\newblock {\em SIAM J. Math. Anal.}, 37(3):889--913, 2005.

\bibitem[M{\'e}t01]{Met}
G.~M{\'e}tivier.
\newblock Stability of multidimensional shocks.
\newblock In {\em Advances in the theory of shock waves}, volume~47 of {\em
  Progr. Nonlinear Differential Equations Appl.}, pages 25--103. Birkh\"{a}user
  Boston, Boston, MA, 2001.

\bibitem[Rod15]{MEDP}
L.~M. Rodrigues.
\newblock Space-modulated stability and averaged dynamics.
\newblock {\em Journ\'ees \'equations aux d\'eriv\'ees partielles}, 2015.
\newblock talk:8.

\bibitem[RZ16]{RZ}
L.~M. Rodrigues and K.~Zumbrun.
\newblock Periodic-coefficient damping estimates, and stability of
  large-amplitude roll waves in inclined thin film flow.
\newblock {\em SIAM J. Math. Anal.}, 48(1):268--280, 2016.

\bibitem[RV98]{RV}
J.-M. Roquejoffre and J.-P. Vila.
\newblock Stability of {ZND} detonation waves in the {M}ajda combustion model.
\newblock {\em Asymptot. Anal.}, 18(3-4):329--348, 1998.

\bibitem[Spe16]{Sp}
J.~Speck.
\newblock {\em Shock formation in small-data solutions to 3{D} quasilinear wave
  equations}, volume 214 of {\em Mathematical Surveys and Monographs}.
\newblock American Mathematical Society, Providence, RI, 2016.

\bibitem[Sze99]{Sz}
A.~Szepessy.
\newblock Dynamics and stability of a weak detonation wave.
\newblock {\em Comm. Math. Phys.}, 202(3):547--569, 1999.

\bibitem[TZ11]{TZ}
B.~Texier and K.~Zumbrun.
\newblock Transition to longitudinal instability of detonation waves is
  generically associated with {H}opf bifurcation to time-periodic galloping
  solutions.
\newblock {\em Comm. Math. Phys.}, 302(1):1--51, 2011.

\bibitem[Wil10]{W}
M.~Williams.
\newblock Heteroclinic orbits with fast transitions: a new construction of
  detonation profiles.
\newblock {\em Indiana Univ. Math. J.}, 59(3):1145--1209, 2010.

\bibitem[YZ20]{YZ}
Z.~Yang and K.~Zumbrun.
\newblock Stability of {H}ydraulic {S}hock {P}rofiles.
\newblock {\em Arch. Ration. Mech. Anal.}, 235(1):195--285, 2020.

\bibitem[Yu99]{YAJOUT}
S.-H. Yu.
\newblock Zero-dissipation limit of solutions with shocks for systems of
  hyperbolic conservation laws.
\newblock {\em Arch. Ration. Mech. Anal.}, 146(4):275--370, 1999.

\bibitem[Zum04]{Z1}
K.~Zumbrun.
\newblock Stability of large-amplitude shock waves of compressible
  {N}avier-{S}tokes equations.
\newblock In {\em Handbook of mathematical fluid dynamics. {V}ol. {III}}, pages
  311--533. North-Holland, Amsterdam, 2004.
\newblock With an appendix by Helge Kristian Jenssen and Gregory Lyng.

\bibitem[Zum07]{Z2}
K.~Zumbrun.
\newblock Planar stability criteria for viscous shock waves of systems with
  real viscosity.
\newblock In {\em Hyperbolic systems of balance laws}, volume 1911 of {\em
  Lecture Notes in Math.}, pages 229--326. Springer, Berlin, 2007.

\bibitem[Zum10]{Z3}
K.~Zumbrun.
\newblock Conditional stability of unstable viscous shock waves in compressible
  gas dynamics and {MHD}.
\newblock {\em Arch. Ration. Mech. Anal.}, 198(3):1031--1056, 2010.

\bibitem[Zum11]{Z5}
K.~Zumbrun.
\newblock Stability of detonation profiles in the {ZND} limit.
\newblock {\em Arch. Ration. Mech. Anal.}, 200(1):141--182, 2011.

\bibitem[Zum12]{Z4}
K.~Zumbrun.
\newblock High-frequency asymptotics and one-dimensional stability of
  {Z}el'dovich--von {N}eumann--{D}\"{o}ring detonations in the small-heat
  release and high-overdrive limits.
\newblock {\em Arch. Ration. Mech. Anal.}, 203(3):701--717, 2012.

\end{thebibliography}

\end{document}